\documentclass[11pt,a4paper]{article}
\usepackage[margin=1in]{geometry}
\usepackage{authblk}
\usepackage{amsmath}
\usepackage{amsfonts}
\usepackage{amsthm}
\usepackage{mathtools}
\usepackage{commath}
\usepackage{graphicx}
\usepackage{caption}
\usepackage{subcaption}
\usepackage[utf8]{inputenc}
\usepackage[english]{babel}
\usepackage{hyperref}
\hypersetup{
    colorlinks=true,
    linkcolor=blue,
    filecolor=magenta,      
    urlcolor=cyan,
    pdftitle={Geometric multiscale analysis via nonstationary subdivision schemes},
}
\usepackage{tikz-cd}
\usetikzlibrary{graphs, decorations.markings}
\graphicspath{{figs/}}
\usepackage{pdflscape}

\newcommand{\N}{\mathbb{N}}

\newtheorem{theorem}{Theorem}[section]

\theoremstyle{definition}
\newtheorem{definition}[theorem]{Definition}
\newtheorem{example}[theorem]{Example}
\newtheorem{remark}[theorem]{Remark}
\newtheorem{corollary}[theorem]{Corollary}

\title{Geometric multiscale analysis via nonstationary subdivision schemes}
\author{Hadar Landau}
\author{Wael Mattar}
\author{Nir Sharon\thanks{Corresponding author: nsharon@tauex.tau.ac.il}}
\affil{School of Mathematical Sciences, Tel Aviv University, Tel Aviv, Israel}
\date{}

\begin{document}

\maketitle

\begin{abstract}
Pyramid transforms are constructive methods for analyzing sequences in a multiscale fashion. Traditionally, these transforms rely on stationary upsampling and downsampling operations. In this paper, we propose employing nonstationary subdivision schemes as upsampling operators that vary with the refinement level. These schemes offer greater flexibility, enabling the development of more expressive multiscale transforms, including geometric multiscale analysis. We establish the fundamental properties of the resulting nonstationary transforms, namely the decay of their detail coefficients and the stability of both the decomposition and the reconstruction, and we demonstrate their effectiveness in capturing and analyzing geometric features. In particular, we apply the framework to assess the consistency of planar samples with a circular shape and to detect localized geometric anomalies, including a study of circles generated by a neural network. The accompanying code is publicly available.
\end{abstract}

\vspace{1em}

\section{Introduction}

Multiscale transforms are essential tools in signal and image processing, enabling a hierarchical mathematical analysis of objects. Customarily, the first scale in the transform corresponds to a coarse representation, and as scales increase, so do the levels of approximation~\cite{daubechies1992ten, mallat1999wavelet}. The pyramid transform uses a refinement or upsampling operator together with a corresponding decimation operator for the construction of a fast multiscale representation of signals in a classical setting that goes back to~\cite{donoho1992interpolating, sweldens1996lifting}. The simplicity of this powerful method opened the door for numerous applications and extensions, such as in geophysics~\cite{li2022simple}, data compression~\cite{lv2018novel}, scientific computing~\cite{gillis2022murphy}, fluid mechanics~\cite{de2021hierarchical}, and geometry~\cite{wallner2020geometric}, just to name a few.

Subdivision schemes are computationally efficient tools for producing smooth objects from discrete sets of points~\cite{dyn1992subdivision}. These schemes are defined by iteratively applying their subdivision operators that refine discrete sets. If a subdivision scheme relies on the same operator at each level of refinement, then it is said to be stationary. Otherwise, it is nonstationary. The subdivision refinements, whether stationary or not, that serve as upsampling operators in the context of pyramid transforms, give rise to a natural connection to multiscale representations~\cite{daubechies1992ten, cohen1996nonstationary}. 

Multiscale transforms, based upon subdivision operators, commonly use interpolating subdivision schemes, i.e., operations that preserve the coarse objects through their refinements. The interpolation simplifies the process of calculating the missing detail coefficients across different scales. In particular, coefficients associated with the interpolating values are systematically zeroed and do not need to be saved or processed at subsequent analysis levels. Recently, the concept of multiscale transforms using subdivision operators has been extended to include noninterpolatory subdivision schemes \cite{dyn2021linear, mattar2023pyramid}. The primary challenge in constructing a noninterpolatory pyramid transform centers on calculating the multiscale details. In particular, given a noninterpolatory subdivision scheme, the corresponding downsampling operators involve applying infinitely supported real-valued sequences. Therefore, care must be taken when realizing and implementing these operators~\cite{mattar2023pyramid, mattar2023pseudo}. 

This paper focuses on pyramid transforms that use nonstationary subdivision schemes. Nonstationary subdivision schemes offer significant advantages over their stationary counterparts, primarily due to their enhanced flexibility and approximation capabilities. For example, unlike stationary schemes that are limited to reproducing algebraic polynomials, nonstationary schemes can reproduce richer function spaces, such as trigonometric, hyperbolic, and exponential polynomials. This property is beneficial for analyzing oscillatory signals, such as speech signals and other band-limited signals in general, or accurately modeling shapes like conic sections and spirals, or for applications involving the solutions of differential equations~\cite{conti2017beyond, conti2011algebraic}. Additionally, the level-dependent nature of nonstationary schemes enables better control over smoothness and geometric features. However, this increased flexibility comes at the cost of greater mathematical and computational complexity, particularly in analyzing convergence and regularity~\cite{dyn1995analysis, dyn2003exponentials}. Despite these challenges, nonstationary schemes provide a powerful tool when stationary schemes are insufficient for the modeling task at hand; for a recent account of the state of the art in nonstationary subdivision, see~\cite{conti2021non}.

The present work differs from the existing subdivision-based pyramid constructions in the level dependence of both components of the transform. The frameworks in~\cite{dyn2021linear,mattar2023pyramid} are built from a fixed subdivision rule, while~\cite{mattar2023pseudo} addresses refinement rules for which an exact or computationally practical even-reverse is not available. In contrast, the present framework employs a sequence of level-dependent refinement operators together with their corresponding level-dependent decimation operators. Consequently, its analysis requires controlling products of nonidentical operators rather than powers of a single stationary operator.

The contributions of this paper are threefold. First, we formulate a general nonstationary pyramid transform based on level-dependent subdivision and decimation pairs. Second, we establish a level-dependent estimate for the detail coefficients and derive stability estimates for both decomposition and reconstruction in terms of the associated transition operators. In particular, summable asymptotic equivalence to a power-bounded stationary scheme yields reconstruction stability uniformly in the number of levels. Third, we exploit the exponential-reproduction capabilities of nonstationary schemes to construct geometry-adapted multiscale diagnostics for perturbations of uniformly sampled circular data and for localized geometric anomalies. More generally, the resulting coefficients measure consistency with the trigonometric or conic model reproduced by the selected scheme.

The structure of the paper is as follows. Section~\ref{sec:back} provides the necessary background and notation, including a concise overview of subdivision schemes and pyramid constructions. In Section~\ref{sec:pyramid}, we introduce the nonstationary framework, present the formal definitions, and establish key properties of nonstationary pyramids. Section~\ref{sec:geom} illustrates the applicability of the proposed framework through geometry-driven pyramids and demonstrates their use in quantifying how closely planar curves resemble circles. Finally, Section~\ref{sec:conclusion} summarizes the main findings and outlines potential directions for future research.

\section{Subdivision schemes and pyramid transforms} \label{sec:back}

This section establishes the foundational definitions and notation required for the developments that follow.

\subsection{Nonstationary subdivision schemes}\label{sec:subdivision_schemes}

Throughout the paper, real-valued sequences will be denoted with bold letters.

A linear univariate subdivision scheme is a sequence of linear subdivision operators $\boldsymbol{\mathcal{S}}=\{\mathcal{S}_{\boldsymbol{\alpha}^{(k)}} \mid   k\in \N \}$ that act on a real-valued bi-infinite sequence $\boldsymbol{c}^{(0)}=\{c_j \in \mathbb{R} \ | \ j \in \mathbb{Z} \}$ as follows
\begin{equation}\label{nonstationary refinement rule}
    \boldsymbol{c}^{(k)} = \mathcal{S}_{\boldsymbol{\alpha}^{(k)}}\boldsymbol{c}^{(k-1)} \quad \text{ where } \quad (\mathcal{S}_{\boldsymbol{\alpha}^{(k)}}\boldsymbol{c}^{(k-1)})_j = \sum_{i \in \mathbb{Z}} \alpha_{j-2i}^{(k)}c_i^{(k-1)}, \quad j \in \mathbb{Z}, \  k \in \N, 
\end{equation}
where $\boldsymbol{\alpha}^{(k)}$ is a finitely supported sequence of real values, called the \emph{mask} used at level $k$ of the subdivision scheme. 
Depending on the parity of the index $j$, we can split the rule~\eqref{nonstationary refinement rule} into the two following rules,
\begin{equation}\label{eqn:two_rules}
\begin{split}
    (\mathcal{S}_{\boldsymbol{\alpha}^{(k)}}\boldsymbol{c}^{(k-1)})_{2j} & =\sum_{i \in \mathbb{Z}} \alpha^{(k)}_{2i}c^{(k-1)}_{j-i}, \\
    (\mathcal{S}_{\boldsymbol{\alpha}^{(k)}}\boldsymbol{c}^{(k-1)})_{2j+1} & =\sum_{i \in \mathbb{Z}} \alpha^{(k)}_{2i+1}c^{(k-1)}_{j-i},
\end{split}
    \qquad j \in \mathbb{Z}, \quad k \in \N.
\end{equation}
A subdivision scheme is said to be \emph{interpolating} if $\alpha^{(k)}_{0}=1$ and $\alpha^{(k)}_{2i}=0$ for all $i\neq0$ and $k\in\N$; equivalently, if $c^{(k)}_{2j}=c^{(k-1)}_{j}$ for all $j\in\mathbb{Z}$, $k \in \N$, and any initial data $\boldsymbol{c}^{(0)}$. 

We define the result after 
$k$ iterations as $ \boldsymbol{c}^{(k)}= \mathcal{S}_{\boldsymbol{\alpha}^{(k)}} \cdots \mathcal{S}_{\boldsymbol{\alpha}^{(1)}} \boldsymbol{c}^{(0)}$. The resulting sequence $\boldsymbol{c}^{(k)}$ is bi-infinite and associated with the grid $2^{-k}\mathbb{Z}$. Consequently, taking the limit as $k$ approaches infinity, we get values $\boldsymbol{\mathcal{S}}^{\infty}\boldsymbol{c}^{(0)}$ that are associated with the dyadic rationals.

A subdivision scheme is termed \emph{convergent} if for any bounded initial data $\boldsymbol{c}^{(0)}$ there exists a continuous function $f$, not identically zero for at least one input, such that
\begin{equation}
    \lim_{k \rightarrow \infty} \sup_{j \in \mathbb{Z}} |f(2^{-k}j)-c^{(k)}_j|=0.
\end{equation}
If such $f$ exists, then it must be unique due to its continuity and the density of the dyadic rationals. In this case, we write $\boldsymbol{\mathcal{S}}^{\infty}\boldsymbol{c}^{(0)}=f$. If, in addition, all the limit functions of the scheme belong to $C^m$ for some $m\in\N\cup\{0\}$, we say that the scheme is $C^m$. If the masks $\boldsymbol{\alpha}^{(k)}$ are the same at all levels, i.e., $\boldsymbol{\alpha}^{(k)}=\boldsymbol{\alpha}$ for all $\ k\in \N$, the scheme is said to be \emph{stationary}; otherwise, it is said to be \emph{nonstationary}.

Nonstationary subdivision schemes offer more flexibility and richer modeling capabilities compared to stationary ones. They can produce a wider variety of limit functions and shapes. For instance, when applied to 2D data, they are capable of generating circles, ellipses, and other conic sections. Moreover, by allowing level-dependent tension parameters, these schemes make it possible to control and fine-tune the shape of the resulting limit function. In the univariate case, they can generate exponential B-splines \cite{dyn1995analysis}, $C^\infty$ functions with compact support such as Rvachev-type functions \cite{dyn2002subdivision}, or B-spline-like functions that achieve higher smoothness relative to the support size \cite{charina2016limits, conti2007totally}.

\begin{remark}
Two additional comments:
    \begin{enumerate}
     \item 
        For brevity, the refinement rule in~\eqref{nonstationary refinement rule} is defined in the linear setting. However, the overall framework readily accommodates non-linear variants as well. In particular, adaptations for manifold-valued data are feasible, in line with existing approaches in the literature (see, e.g.,~\cite{dyn2017manifold}).        
        \item 
        From a practical standpoint, nonstationary subdivision schemes are not significantly more complex than their stationary counterparts, particularly given that only a limited number of refinement levels are typically applied in real-world applications.
    \end{enumerate}
\end{remark}

We recall two examples of nonstationary subdivision schemes, taken from \cite{dyn2003exponentials} and \cite{conti2017beyond}, respectively.

\begin{example}\label{NS example 4-pt}
This scheme is a nonstationary generalization of the interpolatory 4-point Dubuc--Deslauriers scheme, which generates circles. The refinement rules~\eqref{eqn:two_rules} are given by
\begin{equation} \label{NS 4-pt}
\begin{split}
    & c_{2j}^{(k)}=c_j^{(k-1)},\\
    & c_{2j+1}^{(k)}=-w^{(k)}\big(c^{(k-1)}_{j-1}+c^{(k-1)}_{j+2}\big) + \big(1+2\cos(\theta 2^{-k})\big)^2 w^{(k)} \big(c^{(k-1)}_{j}+c^{(k-1)}_{j+1}\big),
\end{split}
\end{equation}
where
\begin{equation*}
    w^{(k)}=\frac{1}{16 \cos^2(\theta 2^{-k-1})\cos(\theta 2^{-k})},
\end{equation*}
and $\theta$ is a scaling factor. Note that the weights of the insertion rule sum to $1$ for every $k$, as follows from the half-angle identity $1+\cos(\theta 2^{-k})=2\cos^2(\theta 2^{-k-1})$. When $\theta=0$, the scheme reduces to the classical 4-point interpolatory scheme based on polynomials~\cite{dyn1991analysis, deslauriers1989symmetric}. Moreover, these insertion rules tend to the $4$-point insertion rules as $k$ approaches infinity. A noteworthy value is $\theta=\frac{2\pi}{N}$: applying the rules to the equally spaced points on the circle with radius $R$ and center at the origin, $\{c_j=R\big(\cos(j\theta),\sin(j\theta)\big)\}_{j=1}^N$, generates a denser set of points on the circle.
\end{example}

\begin{example}\label{NS example cubic B-spline}
This is a nonstationary generalization of the cubic B-spline scheme, which generates exponential splines. The refinement rules~\eqref{eqn:two_rules} are given by
\begin{equation}\label{NS cubic B-spline}
\begin{split}
     & c_{2j}^{(k)}=\frac{1}{2(v^{(k)}+1)^2}c_{j-1}^{(k-1)} + \frac{4(v^{(k)})^2+2}{2(v^{(k)}+1)^2}c_j^{(k-1)} + \frac{1}{2(v^{(k)}+1)^2}c_{j+1}^{(k-1)},\\
     & c_{2j+1}^{(k)}=\frac{2v^{(k)}}{(v^{(k)}+1)^2}c_j^{(k-1)} + \frac{2v^{(k)}}{(v^{(k)}+1)^2}c_{j+1}^{(k-1)} ,
\end{split}
\end{equation}
where the nonstationary parameter $v^{(k)}$ is defined by the recursion
\begin{equation}
    v^{(k)}=\sqrt{\frac{1+v^{(k-1)}}{2}}, \quad k \geq 0, \quad v^{(-1)}>-1.
\end{equation}
Customarily, $v^{(-1)}$ is set to $\cos(\theta)$, where, as before, $\theta$ is a parameter related to the sampling rate; in this case $v^{(k)}=\cos\big(\theta 2^{-(k+1)}\big)$ for all $k\geq0$. Analogously, the hyperbolic initialization $v^{(-1)}=\cosh(\theta)$ yields $v^{(k)}=\cosh\big(\theta 2^{-(k+1)}\big)$. When $\theta=0$, then $v^{(k)}=1$ for all $k$ and the refinement rules coincide with the ones of the cubic B-spline; the same limit is attained as $k$ tends to infinity.
\end{example}

\subsection{Pyramid transforms for data analysis}\label{sec:pyramid_transforms}

The notion of pyramid transforms is constituted by representing a high-resolution sequence of data points as a pyramid composed of a coarse approximation along with layers of details, each associated with a different scale, see for example \cite{donoho1992interpolating, dyn2021linear, grohs2010stability, harten1996multiresolution}.

We begin by reviewing the construction of pyramid transforms using stationary subdivision schemes as upsampling operators, followed by a discussion of various cases determined by the nature of the refinement. Let $\mathcal{S}_{\boldsymbol{\alpha}}$ be an arbitrary subdivision scheme. A high-resolution sequence $\boldsymbol{c}^{(J)}$, $J \in \mathbb{N}$ associated with the fine grid $2^{-J}\mathbb{Z}$, can be recursively decomposed into a coarse (low-resolution) sequence $\boldsymbol{c}^{(0)}$ on the integer grid, along with a collection of details $\boldsymbol{d}^{(\ell)}$, $\ell=1,\dots,J$ that are supported on the grids $2^{-\ell}\mathbb{Z}$, respectively. This representation is called the \emph{pyramid transform} of $\boldsymbol{c}^{(J)}$ via $\mathcal{S}_{\boldsymbol{\alpha}}$ and is done by iterating
\begin{equation} \label{eqn:stationary_pyramid}
    \boldsymbol{c}^{(\ell-1)}=\mathcal{D}_{\boldsymbol{\gamma}}\boldsymbol{c}^{(\ell)}, \quad \boldsymbol{d}^{(\ell)}=\boldsymbol{c}^{(\ell)}-\mathcal{S}_{\boldsymbol{\alpha}}\boldsymbol{c}^{(\ell-1)}, \quad \ell=1,\dots,J,
\end{equation}
where the operator $\mathcal{D}_{\boldsymbol{\gamma}}$ is associated with an absolutely summable bi-infinite sequence $\boldsymbol{\gamma}$ and called the \emph{decimation} operator. It was first introduced in~\cite{dyn2021linear} and is explicitly given by 
\begin{equation}\label{decimation operator}
    (\mathcal{D}_{\boldsymbol{\gamma}}\boldsymbol{c})_j = \sum_{i \in \mathbb{Z}}\gamma_{j-i}c_{2i}, \quad j\in \mathbb{Z},
\end{equation}
for any real-valued sequence $\boldsymbol{c}$. The operator $\mathcal{D}_{\boldsymbol{\gamma}}(\boldsymbol{c})$ can be realized as a convolution between $\boldsymbol{\gamma}$ and the even-indexed elements of $\boldsymbol{c}$. Decimation operators serve as downscaling tools, reducing the number of data points when applied to a sequence.

When the recursive process~\eqref{eqn:stationary_pyramid} acts on a sequence $\boldsymbol{c}^{(J)}$, it produces a pyramid of sequences denoted by $\{\boldsymbol{c}^{(0)}; \boldsymbol{d}^{(1)}, \dots, \boldsymbol{d}^{(J)}\}$. Conversely, the \emph{inverse} pyramid transform, usually termed as the \emph{synthesis}, takes the pyramid representation and reconstructs the original sequence $\boldsymbol{c}^{(J)}$ iteratively by
\begin{equation} \label{synthesis}
    \boldsymbol{c}^{(\ell)}= \mathcal{S}_{\boldsymbol{\alpha}}\boldsymbol{c}^{(\ell-1)}+\boldsymbol{d}^{(\ell)}, \quad \ell=1,\dots,J.
\end{equation}

The primary purpose of decimation operators is to ensure that the even detail coefficients produced by~\eqref{eqn:stationary_pyramid} are zero. In particular, these operators were mainly designed to satisfy $d^{(\ell)}_{2j}=0$ for all $\ell=1,\dots,J$ and $j\in\mathbb{Z}$~\cite{dyn2021linear}. This requirement is the backbone of pyramid transforms that are based on refinement operators. In particular, this property is mandatory for an application like \emph{data compression} since one can discard ``half'' of the detail coefficients when storing the pyramid representation.

The necessity of having zero detail coefficients at even indices across all levels implies the following relation between the mask $\boldsymbol{\alpha}$ and the sequence $\boldsymbol{\gamma}$ involved in~\eqref{eqn:stationary_pyramid}. With a subdivision scheme mask $\boldsymbol{\alpha}$ at hand, the even detail coefficients of~\eqref{eqn:stationary_pyramid} vanish if the sequence $\boldsymbol{\gamma}$ satisfies the convolutional equation
\begin{equation} \label{convolutional equation for finding gamma}
    \boldsymbol{\gamma}*(\boldsymbol{\alpha}\downarrow2)=\boldsymbol{\delta},
\end{equation}
where $\boldsymbol{\delta}$ is the Kronecker delta sequence ($\delta_0=1$ and $\delta_i=0$ for $i \neq 0$), and $\downarrow 2$ is the simple downsampling operator that returns the even elements of a sequence. Namely, $(\boldsymbol{c}\downarrow2)_{j}=c_{2j}$ for all $j\in\mathbb{Z}$.

If the subdivision scheme $\mathcal{S}_{\boldsymbol{\alpha}}$ is interpolating, then $\boldsymbol{\alpha}\downarrow2=\boldsymbol{\delta}$ and hence $\boldsymbol{\gamma}=\boldsymbol{\delta}$. In this case, the decimation operator $\mathcal{D}_{\boldsymbol{\delta}}$ reduces to the rudimentary $\downarrow2$ operation. Otherwise, the well-known Wiener's lemma, see~\cite[Chapter 5]{forster2010four}, provides a sufficient condition for the existence of an absolutely summable solution of~\eqref{convolutional equation for finding gamma}: if $\boldsymbol{\alpha}\downarrow2\in\ell^1(\mathbb{Z})$ and its symbol does not vanish on the unit circle, that is, $\sum_{i\in\mathbb{Z}}(\boldsymbol{\alpha}\downarrow2)_i z^i\neq0$ for all $|z|=1$, then $\boldsymbol{\alpha}\downarrow2$ is invertible with respect to convolution, and~\eqref{convolutional equation for finding gamma} admits a unique solution $\boldsymbol{\gamma}\in\ell^1(\mathbb{Z})$. Moreover, the recent study~\cite{mattar2023pseudo} has introduced a method called \emph{pseudo-reversing} for approximating $\boldsymbol{\gamma}$ in case it is not absolutely summable.

In general, if a solution to~\eqref{convolutional equation for finding gamma} exists, then it is infinitely supported~\cite{grochenig2010wiener}, and its elements must sum to $1$ provided that $\sum_{i\in\mathbb{Z}}\alpha_{2i}=1$. In this case, following~\cite{dyn2021linear}, we say that $\mathcal{D}_{\boldsymbol{\gamma}}$ is the \emph{even-reverse} of $\mathcal{S}_{\boldsymbol{\alpha}}$. Furthermore, there exist constants $C>0$ and $\lambda\in(0,1)$ depending on the size of $\boldsymbol{\alpha}\downarrow2$ such that
\begin{equation} \label{gamma's decay}
    |\gamma_j| \leq C\lambda^{|j|}, \quad j \in \mathbb{Z}.
\end{equation}
In fact, these constants depend on the \emph{reversibility condition number}~\cite{mattar2023pseudo} that takes values greater or equal to $1$, including $\infty$ for non-reversible masks. For compactness we avoid the explicit formulae of the constants. However, similar to the notion of condition number in linear algebra, the lower the reversibility condition number is, the better the decay of $\boldsymbol{\gamma}$ becomes. Overall, the geometric bound~\eqref{gamma's decay} is essential for the computation of $\mathcal{D}_{\boldsymbol{\gamma}}$ as we will discuss next.

In theory, the application of the operator~\eqref{decimation operator} involves an infinite amount of computations due to the infinite support of the sequence $\boldsymbol{\gamma}$. However, in practice, we approximate the operator by truncating its sequence $\boldsymbol{\gamma}$ as shown in~\cite{mattar2023pyramid}. In particular, the decimation operator $\mathcal{D}_{\boldsymbol{\gamma}}$ can be approximated via a proper truncation of its associated sequence. Given a \emph{truncation} parameter $\varepsilon>0$, we only consider the elements of $\boldsymbol{\gamma}$ that are greater (in absolute value) than $\varepsilon$. The other elements are therefore negligible and treated as zeros. The bound~\eqref{gamma's decay} guarantees that the support of the resulting truncation is finite for any positive $\varepsilon$. Consequently, the truncated sequence is normalized so that its elements sum to $1$. The normalization is done by element-wise division of the truncated sequence by its total sum. This type of normalization is essential for multiscaling as it directly affects the significance of the even detail coefficients, as shown in~\cite{mattar2023pyramid}.

Once truncated and normalized, the resulting approximating sequence then replaces $\boldsymbol{\gamma}$ in $\mathcal{D}_{\boldsymbol{\gamma}}$, therefore making the transform~\eqref{eqn:stationary_pyramid} implementable with controllable errors depending on the truncation parameter $\varepsilon$. As shown in~\cite{mattar2023pyramid}, the resulting multiscale transform enjoys two main properties, the decay of the detail coefficients and the stability of the inverse transform. In the next section, we introduce and study the nonstationary version of the transform~\eqref{eqn:stationary_pyramid}.


\section{Pyramid transforms based on nonstationary schemes}\label{sec:pyramid}

This section is devoted to pyramid transforms based on nonstationary subdivision schemes, that is, transforms whose upsampling operators vary with the level. The resulting construction extends the classical stationary framework of Subsection~\ref{sec:pyramid_transforms}. To illustrate the construction and potential of this approach, two representative examples are presented, followed by the fundamental properties of the resulting transforms.

Pyramid transforms can be implemented in a level-dependent manner by utilizing nonstationary subdivision schemes~\eqref{nonstationary refinement rule}. In particular, multiscaling a real-valued sequence $\boldsymbol{c}^{(J)}$ with a nonstationary family of subdivision schemes $\{\mathcal{S}_{\boldsymbol{\alpha}^{(\ell)}} \mid  \ell=1,\dots,J\}$ is defined by iterating
\begin{equation} \label{NS multiscale transform}
    \boldsymbol{c}^{(\ell-1)}=\mathcal{D}_{\boldsymbol{\gamma}^{(\ell)}}\boldsymbol{c}^{(\ell)}, \quad \boldsymbol{d}^{(\ell)}=\boldsymbol{c}^{(\ell)}-\mathcal{S}_{\boldsymbol{\alpha}^{(\ell)}}\boldsymbol{c}^{(\ell-1)}, \quad \ell=1,\dots,J,
\end{equation}
where $\mathcal{D}_{\boldsymbol{\gamma}^{(\ell)}}$ is the decimation operator~\eqref{decimation operator} corresponding to $\mathcal{S}_{\boldsymbol{\alpha}^{(\ell)}}$. That is, the pair $\boldsymbol{\gamma}^{(\ell)}$ and $\boldsymbol{\alpha}^{(\ell)}$ solve~\eqref{convolutional equation for finding gamma}. Conversely, similar to~\eqref{synthesis}, the synthesis is obtained by iterating
\begin{equation} \label{NS synthesis}
    \boldsymbol{c}^{(\ell)}= \mathcal{S}_{\boldsymbol{\alpha}^{(\ell)}}\boldsymbol{c}^{(\ell-1)}+\boldsymbol{d}^{(\ell)}, \quad \ell=1,\dots,J.
\end{equation}

We illustrate the nonstationary pyramid in Figure~\ref{fig:epyramid}.
\begin{figure}[htbp] 
    \begin{center}
\begin{tikzcd}[column sep=large, row sep=large]
    \boldsymbol{c}^{(J)} \arrow[d, "\mathcal{D}_{\boldsymbol{\gamma}^{(J)}}"'] \arrow[dr, "-"] & \\
    \boldsymbol{c}^{(J-1)} \arrow[d, "\mathcal{D}_{\boldsymbol{\gamma}^{(J-1)}}"'] \arrow[dr, "-"] \arrow[r, "\mathcal{S}_{\boldsymbol{\alpha}^{(J)}}"] & \boldsymbol{d}^{(J)} \\
    \boldsymbol{c}^{(J-2)} \arrow[d, dashed] \arrow[r, "\mathcal{S}_{\boldsymbol{\alpha}^{(J-1)}}"] & \boldsymbol{d}^{(J-1)} \\
    \boldsymbol{c}^{(1)} \arrow[d, "\mathcal{D}_{\boldsymbol{\gamma}^{(1)}}"'] \arrow[dr, "-"] & \\
    \boldsymbol{c}^{(0)} \arrow[r, "\mathcal{S}_{\boldsymbol{\alpha}^{(1)}}"] & \boldsymbol{d}^{(1)}
\end{tikzcd}
\qquad\qquad
\begin{tikzcd}[column sep=large, row sep=large]
    & \boldsymbol{c}^{(J)} \\
    \boldsymbol{d}^{(J)} \arrow[ur,"\quad+" '] 
    & \boldsymbol{c}^{(J-1)} \arrow[u, "\mathcal{S}_{\boldsymbol{\alpha}^{(J)}}" '] \\
    & \boldsymbol{c}^{(2)} \arrow[u, dashed] \\
    \boldsymbol{d}^{(2)} \arrow[ur,"\quad+" '] 
        & \boldsymbol{c}^{(1)} \arrow[u, "\mathcal{S}_{\boldsymbol{\alpha}^{(2)}}" '] \\
    \boldsymbol{d}^{(1)} \arrow[ur,"\quad+" '] 
        & \boldsymbol{c}^{(0)} \arrow[u, "\mathcal{S}_{\boldsymbol{\alpha}^{(1)}}" ']
\end{tikzcd}
\caption{Schematic depiction of the nonstationary pyramid transform. On the left is the analysis~\eqref{NS multiscale transform}, and on the right is the synthesis~\eqref{NS synthesis}.}
 \label{fig:epyramid}
\end{center}
\end{figure}

Although we use a different decimation operator at each decomposition iteration, we truncate the sequences $\boldsymbol{\gamma}^{(\ell)}$, $\ell=1,\dots, J$ with a uniform truncation threshold $\varepsilon$. The truncated sequences are then normalized as described in Subsection~\ref{sec:pyramid_transforms}. The resulting sequences then replace $\boldsymbol{\gamma}^{(\ell)}$ in~\eqref{NS multiscale transform}. We proceed by discussing the nonstationary pyramid transforms corresponding to Examples~\ref{NS example 4-pt} and~\ref{NS example cubic B-spline}.

\begin{example}
The subdivision scheme presented in \eqref{NS 4-pt} is interpolatory, and therefore its even refinement rule is stationary. Consequently, the downscaling operators $\mathcal{D}_{\boldsymbol{\gamma}^{(\ell)}}$, $\ell=1,\dots,J$ reduce to the simple downsampling operator $\downarrow2$. Hence, in the nonstationary multiscaling~\eqref{NS multiscale transform} the computation of $\boldsymbol{c}^{(\ell-1)}$ simplifies to $\boldsymbol{c}^{(\ell-1)}=\boldsymbol{c}^{(\ell)}\downarrow 2$. 
\end{example}

\begin{figure}[htbp]
    \centering
    \begin{subfigure}{0.24\textwidth}
         \centering         \includegraphics[width=\textwidth]{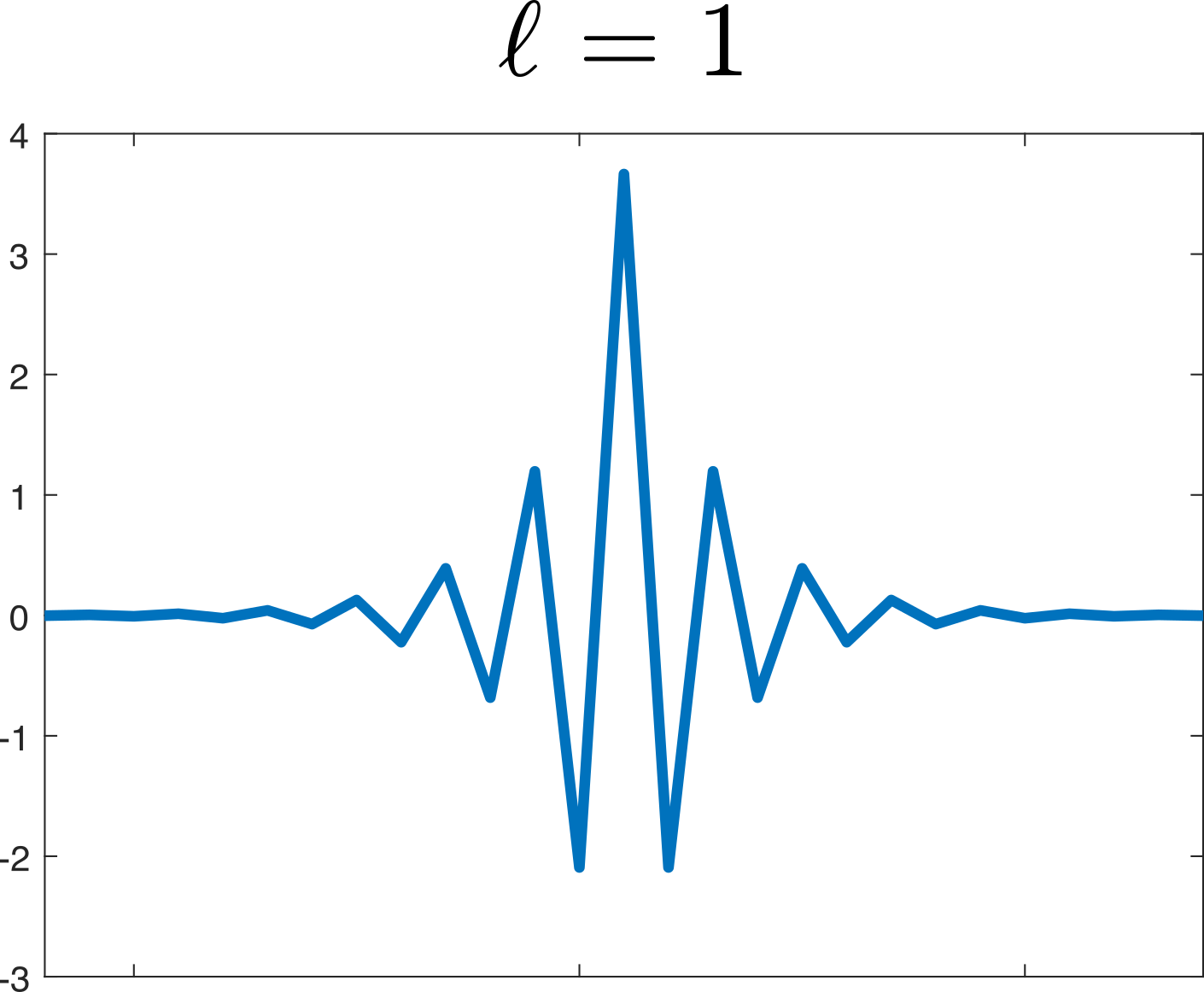}
    \end{subfigure} 
    \begin{subfigure}{0.24\textwidth}
         \centering         \includegraphics[width=\textwidth]{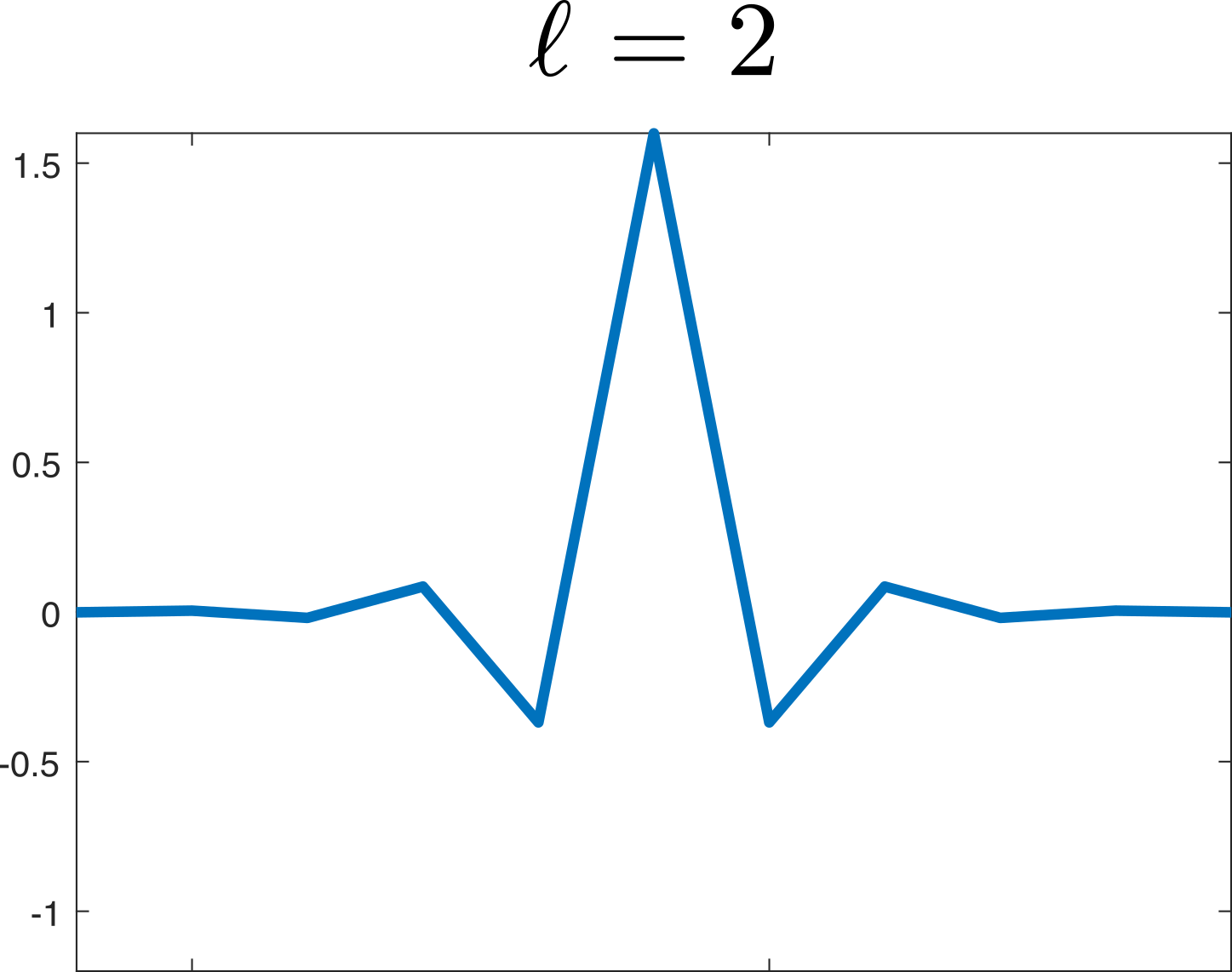}
     \end{subfigure}
    \begin{subfigure}{0.24\textwidth}
         \centering         \includegraphics[width=\textwidth]{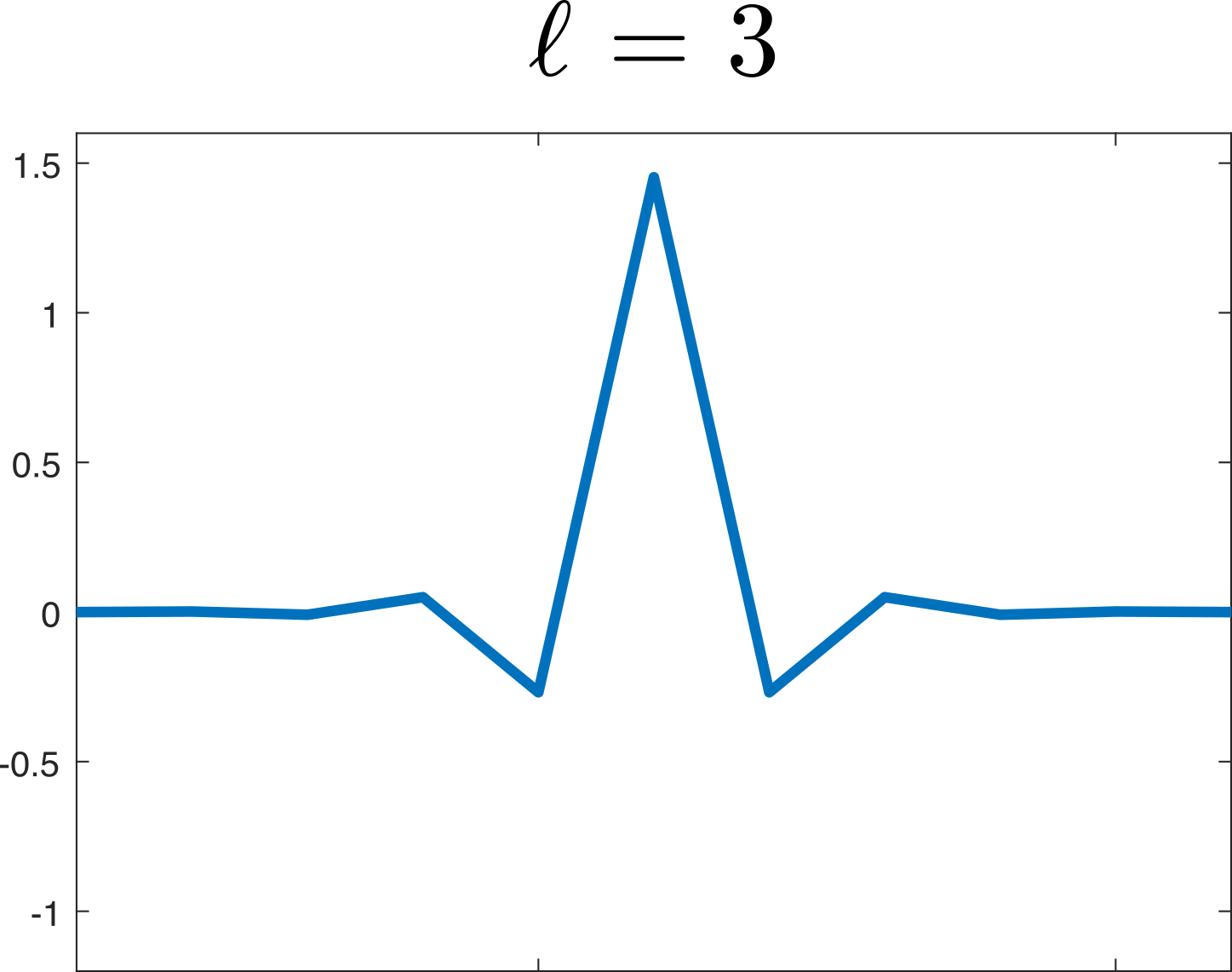}
     \end{subfigure}
    \begin{subfigure}{0.24\textwidth}
         \centering         \includegraphics[width=\textwidth]{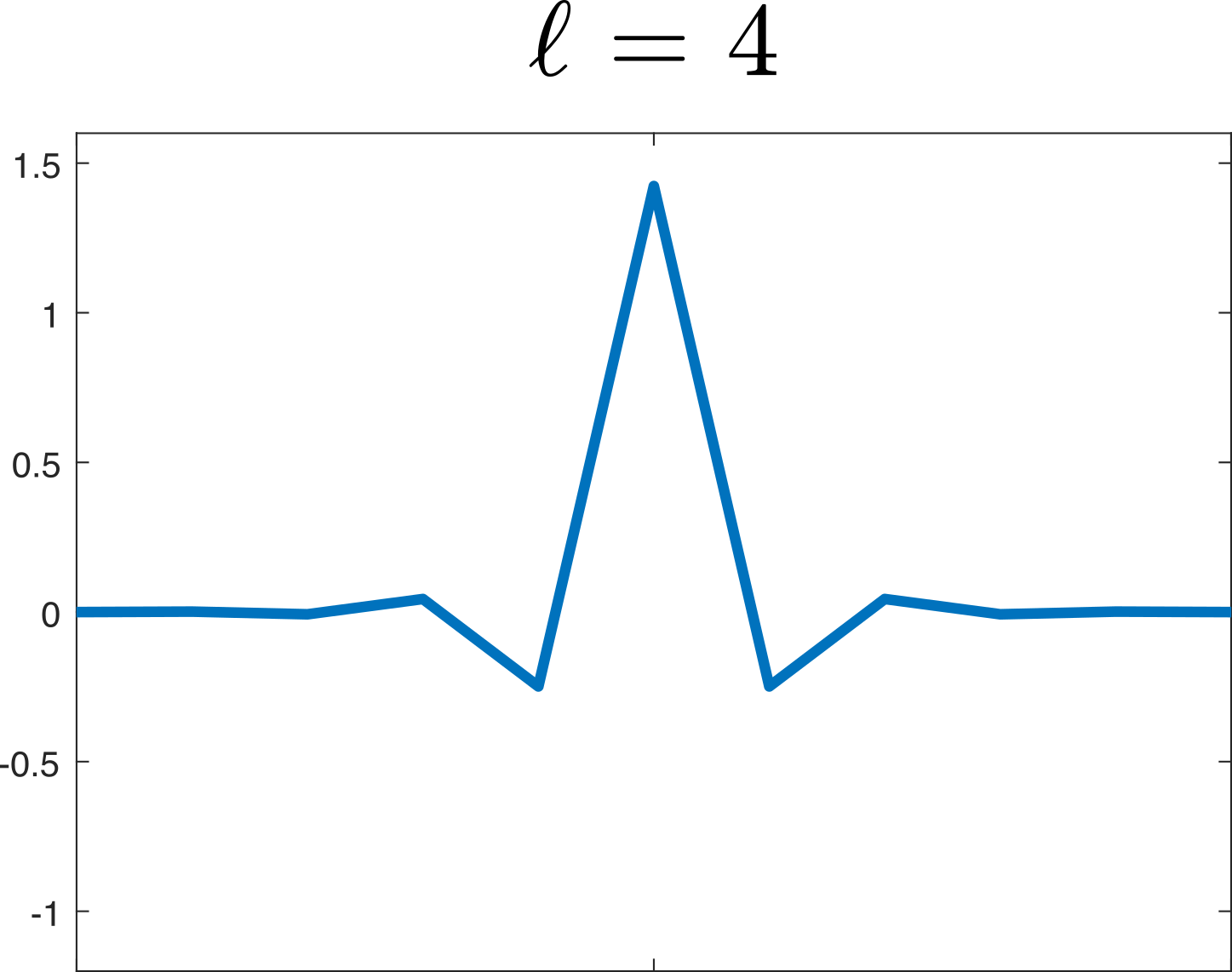}
     \end{subfigure}
     \caption{The coefficients of the first four decimation operators corresponding to the nonstationary subdivision scheme defined in~\eqref{NS cubic B-spline}. For each level $\ell$, the decimation operator is obtained by solving~\eqref{convolutional equation for finding gamma}, followed by truncation at a fixed threshold $\varepsilon = 10^{-15}$ and subsequent normalization. A notable difference is observed between the first and second operators, while the remaining ones display similar profiles and progressively tend toward the (stationary) cubic B-spline decimation operator.}
     \label{fig:NS decimation}
\end{figure}

\begin{example}
The subdivision scheme given in~\eqref{NS cubic B-spline} is noninterpolatory. As a result, the corresponding decimation operators $\mathcal{D}_{\boldsymbol{\gamma}^{(\ell)}}$, for $\ell = 1, \dots, J$, are associated with sequences with infinite support. Notably, as the level $\ell$ increases, the operator $\mathcal{D}_{\boldsymbol{\gamma}^{(\ell)}}$ asymptotically approaches the cubic B-spline decimation operator. This behavior arises from the continuity of the solution of~\eqref{convolutional equation for finding gamma} with respect to the refinement mask of the subdivision scheme. The coefficients of the resulting decimation operators at the first four levels are illustrated in \figref{fig:NS decimation}.
\end{example}

It turns out that stationary and nonstationary transforms share similar properties. Let us now present the decay of the details and the stability theorems for the nonstationary case. The following results and their proofs are analogous to those shown in~\cite{mattar2023pyramid} and \cite{dyn2021linear}.

\begin{theorem}[Decay of the detail coefficients]\label{thm. details decay}
Let $f\colon \mathbb{R} \rightarrow \mathbb{R}$ be a bounded, continuously differentiable function with a bounded derivative, and let $\boldsymbol{c}^{(J)}$, $J \in \mathbb{N}$, be its samples over the dyadic grid $2^{-J}\mathbb{Z}$, that is, $\boldsymbol{c}^{(J)}=f|_{2^{-J}\mathbb{Z}}$. Let $\{\mathcal{S}_{\boldsymbol{\alpha}^{(\ell)}} \mid \ell=1,\dots,J\}$ be a family of nonstationary subdivision schemes satisfying $\sum_{i}\alpha^{(\ell)}_{2i}=\sum_{i}\alpha^{(\ell)}_{2i+1}=1$, and let $\{\mathcal{D}_{\boldsymbol{\gamma}^{(\ell)}} \mid \ell=1,\dots,J\}$ be its corresponding family of decimation operators, associated with sequences $\boldsymbol{\gamma}^{(\ell)}$ satisfying $\sum_i \gamma^{(\ell)}_i = 1$ and having finite first moments, that is, $\sum_{i\in\mathbb{Z}}|i|\,|\gamma^{(\ell)}_i|<\infty$. Then the detail coefficients obtained by the multiscale transform~\eqref{NS multiscale transform} satisfy
\begin{equation}\label{details bound}
    \|\boldsymbol{d}^{(\ell)}\|_\infty \leq \bigg(K_{\boldsymbol{\alpha}^{(\ell)},\boldsymbol{\gamma}^{(\ell)}} \cdot \|f^\prime\|_\infty \cdot \prod_{m=\ell+1}^J \|\boldsymbol{\gamma}^{(m)}\|_1 \bigg) \; 2^{-\ell}, \quad \ell=1,\dots,J,
\end{equation}
where $K_{\boldsymbol{\alpha}^{(\ell)},\boldsymbol{\gamma}^{(\ell)}}=K_{\boldsymbol{\gamma}^{(\ell)}}\|\boldsymbol{\alpha}^{(\ell)}\|_1+ K_{\boldsymbol{\alpha}^{(\ell)}}\|\boldsymbol{\gamma}^{(\ell)}\|_1$, with $K_{\boldsymbol{\gamma}^{(\ell)}}=2\sum_{i\in\mathbb{Z}} |\gamma^{(\ell)}_i|\,|i|$ and $K_{\boldsymbol{\alpha}^{(\ell)}}=\sum_{i\in\mathbb{Z}} |\alpha^{(\ell)}_i|\,|i|$, and with the convention that the empty product, occurring for $\ell=J$, equals $1$.
\end{theorem}

\begin{proof}
    Throughout, for a bounded sequence $\boldsymbol{c}$ we write $\Delta\boldsymbol{c}\mathrel{:=}\sup_{k\in\mathbb{Z}}|c_{k+1}-c_k|$ for its maximal first difference; note that $\Delta$ is subadditive and satisfies $\Delta(\boldsymbol{\gamma}*\boldsymbol{c})\leq\|\boldsymbol{\gamma}\|_1\,\Delta\boldsymbol{c}$ for any absolutely summable sequence $\boldsymbol{\gamma}$, since the difference operator commutes with convolution. We first show that $\|\boldsymbol{d}^{(\ell)}\|_\infty$ is proportional to $\Delta\boldsymbol{c}^{(\ell)}$ for all levels $\ell=1,\dots,J$, namely
    \begin{equation}\label{eqn:details_Delta}
        \|\boldsymbol{d}^{(\ell)}\|_\infty \leq K_{\boldsymbol{\alpha}^{(\ell)},\boldsymbol{\gamma}^{(\ell)}} \Delta \boldsymbol{c}^{(\ell)}.
    \end{equation}
    To this end, we begin by evaluating a general term in $\boldsymbol{d}^{(\ell)}$. Observe
    \begin{equation*}
    \begin{split}
        d^{(\ell)}_k & =c^{(\ell)}_k-\sum_i \alpha^{(\ell)}_{k-2i}c^{(\ell-1)}_i=\sum_i \alpha^{(\ell)}_{k-2i}(c^{(\ell)}_k-c^{(\ell-1)}_i)\\
        & =\sum_i \alpha^{(\ell)}_{k-2i}\bigg(c^{(\ell)}_k-\sum_n \gamma^{(\ell)}_{i-n} c^{(\ell)}_{2n}\bigg)\\
        & =\sum_i \alpha^{(\ell)}_{k-2i}\bigg(\sum_n \gamma^{(\ell)}_{i-n}c^{(\ell)}_k-\sum_n \gamma^{(\ell)}_{i-n} c^{(\ell)}_{2n}\bigg)\\
        & =\sum_i \alpha^{(\ell)}_{k-2i}\bigg(\sum_n \gamma^{(\ell)}_{i-n}(c^{(\ell)}_k-c^{(\ell)}_{2n})\bigg).
    \end{split}
    \end{equation*}
    Consequently,
    \begin{equation*}
    \begin{split}
        |d^{(\ell)}_k| & \leq \sum_i|\alpha^{(\ell)}_{k-2i}|\bigg(\sum_n |\gamma^{(\ell)}_{i-n}|\cdot |c^{(\ell)}_k-c^{(\ell)}_{2n}|\bigg)\\ & \leq \sum_i|\alpha^{(\ell)}_{k-2i}|\bigg(\sum_n |\gamma^{(\ell)}_{i-n}|\cdot |2n-k| \cdot \Delta\boldsymbol{c}^{(\ell)}\bigg)\\
        & \leq \sum_i|\alpha^{(\ell)}_{k-2i}|\bigg(\sum_n |\gamma^{(\ell)}_{i-n}|\cdot \big(|2n-2i|+|2i-k|\big) \bigg) \cdot \Delta\boldsymbol{c}^{(\ell)}\\
        & \leq \sum_i|\alpha^{(\ell)}_{k-2i}|\bigg(K_{\boldsymbol{\gamma}^{(\ell)}}+|2i-k| \|\boldsymbol{\gamma}^{(\ell)}\|_1\bigg) \cdot \Delta\boldsymbol{c}^{(\ell)}\\
        & \leq \Big(K_{\boldsymbol{\gamma}^{(\ell)}}\|\boldsymbol{\alpha}^{(\ell)}\|_1+K_{\boldsymbol{\alpha}^{(\ell)}} \|\boldsymbol{\gamma}^{(\ell)}\|_1\Big) \cdot \Delta\boldsymbol{c}^{(\ell)},\\
    \end{split}
    \end{equation*}
    where the last inequality follows from $\sum_i|\alpha^{(\ell)}_{k-2i}|\leq\|\boldsymbol{\alpha}^{(\ell)}\|_1$ and $\sum_i|\alpha^{(\ell)}_{k-2i}|\,|k-2i|\leq K_{\boldsymbol{\alpha}^{(\ell)}}$, since for a fixed $k$ the summation runs over a single parity class of the mask.

    Applying the supremum over $k\in\mathbb{Z}$ to the latter inequality yields~\eqref{eqn:details_Delta}. We now evaluate the term $\Delta\boldsymbol{c}^{(\ell)}$ for an arbitrary $\ell$ based on the assumption that $\boldsymbol{c}^{(J)}$ is sampled from $f$ on the equispaced grid $2^{-J}\mathbb{Z}$. Since $f$ is bounded and differentiable, on the finest level, we have that
    \begin{equation}\label{eqn:Delta_c^J}
        \Delta\boldsymbol{c}^{(J)}\leq 2^{-J}\|f^\prime\|_{\infty}.
    \end{equation}
    This inequality is guaranteed by the Mean Value Theorem. In particular, for any $k\in\mathbb{Z}$ there exists a point $x_k$ in the open interval connecting the two consecutive grid points $k2^{-J}$ and $(k+1)2^{-J}$, such that $|c^{(J)}_{k+1}-c^{(J)}_{k}|\leq 2^{-J}|f^\prime (x_k)|$, and hence~\eqref{eqn:Delta_c^J} is obtained by taking the supremum over the index $k$. On coarser levels, however, since $\Delta$ commutes with the convolution operator, we get the recursive relation

    \begin{equation}\label{eqn:Delta_recursion}
        \Delta\boldsymbol{c}^{(\ell-1)} = \Delta\big(\boldsymbol{\gamma}^{(\ell)}*(\boldsymbol{c}^{(\ell)}\downarrow 2)\big)\leq \|\boldsymbol{\gamma}^{(\ell)}\|_1\cdot\Delta(\boldsymbol{c}^{(\ell)}\downarrow 2)\leq 2\|\boldsymbol{\gamma}^{(\ell)}\|_1\cdot\Delta\boldsymbol{c}^{(\ell)},
    \end{equation}
    which holds for $\ell=1,\dots,J$. Finally, iterative calculations of~\eqref{eqn:Delta_recursion} involving~\eqref{eqn:Delta_c^J} and~\eqref{eqn:details_Delta} yield~\eqref{details bound}.
\end{proof}

Theorem~\ref{thm. details decay} provides the level-dependent estimate \eqref{details bound}. To obtain an \(O(2^{-\ell})\) bound uniformly in the number of levels, it is sufficient to assume
$$ \sup_{\ell\ge1}K_{\alpha^{(\ell)},\gamma^{(\ell)}}<\infty \quad\text{and}\quad \sup_{J\ge1}\max_{1\le\ell\le J} \prod_{m=\ell+1}^{J}\|\gamma^{(m)}\|_1<\infty. $$
A sufficient condition for the second bound is
$$ \sum_{m=1}^{\infty}\bigl(\|\gamma^{(m)}\|_1-1\bigr)<\infty. $$
Under these conditions, \eqref{details bound} yields \(\|d^{(\ell)}\|_\infty=O(2^{-\ell})\), uniformly in \(J\). In the stationary case, the constant becomes level-independent, and the product reduces to \(\|\gamma\|_1^{J-\ell}\), recovering Corollary 3.5 in~\cite{mattar2023pyramid}.

\begin{remark} \label{rem:hypotheses}
Three comments regarding the hypotheses of Theorem~\ref{thm. details decay} are in order. First, the proof uses only the normalization $\sum_i \gamma^{(\ell)}_i=1$ together with the sum rules of the masks; the even-reverse relation~\eqref{convolutional equation for finding gamma} is not required. Hence, the bound~\eqref{details bound} is valid for any normalized decimation operators. Second, the sum rules $\sum_{i}\alpha^{(\ell)}_{2i}=\sum_{i}\alpha^{(\ell)}_{2i+1}=1$ hold exactly for schemes reproducing linear polynomials, such as the interpolatory scheme of Example~\ref{NS example 4-pt} and the scheme of Example~\ref{example:geometric_scheme} below. However, exponential-generating schemes may satisfy them only asymptotically. For instance, for the scheme of Example~\ref{NS example cubic B-spline}, the parity sums of the mask are $1\pm(1-v^{(\ell)})^2(1+v^{(\ell)})^{-2}$, deviating from $1$ by $O(16^{-\ell})$ when $v^{(\ell)}=\cos\big(\theta 2^{-(\ell+1)}\big)$. A routine modification of the proof shows that, in this case,~\eqref{details bound} holds with an additional additive term proportional to this deviation times $\|\boldsymbol{c}^{(\ell)}\|_\infty$. Since $\|\boldsymbol{c}^{(\ell)}\|_\infty\leq\prod_{m=\ell+1}^{J}\|\boldsymbol{\gamma}^{(m)}\|_1\,\|f\|_\infty$, the norms $\|\boldsymbol{c}^{(\ell)}\|_\infty$ are bounded uniformly in $\ell$ whenever the products $\prod_{m=\ell+1}^{J}\|\boldsymbol{\gamma}^{(m)}\|_1$ are, as in the near-interpolating regime discussed above. Under this uniform bound, the additional term is of order $16^{-\ell}$ and decays faster than the leading $2^{-\ell}$ factor, leaving the overall decay unaffected. Third, the theorem concerns the exact solutions $\boldsymbol{\gamma}^{(\ell)}$ of~\eqref{convolutional equation for finding gamma}; in practice these sequences are truncated and renormalized, and the effect of the truncation on the constants is controlled by the truncation parameter $\varepsilon$, analogously to the stationary analysis in~\cite{mattar2023pyramid}.
\end{remark}

Let us now show that the inverse multiscale transform~\eqref{NS synthesis} is stable. The stability estimate is naturally expressed in terms of the transition operators
\begin{equation}\label{eqn:transition_operators}
    \mathcal{P}_{q,p}=\mathcal{S}_{\boldsymbol{\alpha}^{(q)}}\mathcal{S}_{\boldsymbol{\alpha}^{(q-1)}}\cdots\mathcal{S}_{\boldsymbol{\alpha}^{(p)}}, \quad 1\leq p\leq q\leq J, \qquad \mathcal{P}_{p-1,p}=\mathcal{I},
\end{equation}
where $\mathcal{I}$ denotes the identity operator; that is, the empty product is understood as the identity. We further denote
\begin{equation}\label{eqn:LJ}
    L_J=\max\Big\{1,\; \max_{1\leq p\leq q\leq J}\|\mathcal{P}_{q,p}\|_\infty\Big\},
\end{equation}
so that $L_J$ bounds the norms of all the transition operators, including the empty ones.

\begin{theorem}[Stability of the reconstruction]
\label{thm. stability of reconstruction}
Let $\{\boldsymbol{c}^{(0)};\boldsymbol{d}^{(1)},\dots,\boldsymbol{d}^{(J)}\}$ and $\{\widetilde{\boldsymbol{c}}^{(0)};\widetilde{\boldsymbol{d}}^{(1)},\dots,\widetilde{\boldsymbol{d}}^{(J)}\}$ be two arbitrary sets of pyramid data, for instance, obtained by~\eqref{NS multiscale transform}, possibly followed by a manipulation of the coefficients. Let $\{\mathcal{S}_{\boldsymbol{\alpha}^{(\ell)}} \mid \ell=1,\dots,J\}$ be a family of nonstationary subdivision schemes with transition operators~\eqref{eqn:transition_operators}. Then the sequences $\boldsymbol{c}^{(J)}$ and $\widetilde{\boldsymbol{c}}^{(J)}$, reconstructed via~\eqref{NS synthesis} from the two pyramids, respectively, obey
\begin{equation}\label{eqn:reconstruction_stability}
    \|\boldsymbol{c}^{(J)}-\widetilde{\boldsymbol{c}}^{(J)}\|_\infty \leq L_J \bigg(\|\boldsymbol{c}^{(0)}-\widetilde{\boldsymbol{c}}^{(0)}\|_\infty+ \sum_{\ell=1}^J \|\boldsymbol{d}^{(\ell)}-\widetilde{\boldsymbol{d}}^{(\ell)}\|_\infty \bigg),
\end{equation}
with $L_J$ of~\eqref{eqn:LJ}. In particular, if
\begin{equation}\label{eqn:uniform_L}
    L=\sup_{J\geq1}L_J<\infty,
\end{equation}
then~\eqref{eqn:reconstruction_stability} holds with $L_J$ replaced by $L$, uniformly in the number of levels $J$.
\end{theorem}

\begin{proof}
    By the linearity of the subdivision operators~\eqref{nonstationary refinement rule}, iterating the reconstruction formula~\eqref{NS synthesis} yields the identity
    \begin{equation*}
        \boldsymbol{c}^{(J)}-\widetilde{\boldsymbol{c}}^{(J)} = \mathcal{P}_{J,1}\big(\boldsymbol{c}^{(0)}-\widetilde{\boldsymbol{c}}^{(0)}\big) + \sum_{\ell=1}^J \mathcal{P}_{J,\ell+1}\big(\boldsymbol{d}^{(\ell)}-\widetilde{\boldsymbol{d}}^{(\ell)}\big),
    \end{equation*}
    where the term $\ell=J$ involves the empty product $\mathcal{P}_{J,J+1}=\mathcal{I}$. Applying the triangle inequality and bounding the norm of each transition operator by $L_J$ of~\eqref{eqn:LJ} gives~\eqref{eqn:reconstruction_stability}. The uniform statement follows immediately from the definition of $L$ in~\eqref{eqn:uniform_L}.
\end{proof}

We note that the bound~\eqref{eqn:reconstruction_stability} involves norms of products of subdivision operators rather than products of their norms. Since the operator norm is submultiplicative, $L_J\leq\max\{1,M\}^J$ with $M=\max_{\ell=1,\dots,J}\|\mathcal{S}_{\boldsymbol{\alpha}^{(\ell)}}\|_\infty$; however, $L_J$ can be significantly smaller and, in contrast to the latter exponential bound, it may remain bounded as $J$ grows. The operator norm of a subdivision operator is computed explicitly from its mask via
\begin{equation}\label{eqn:subdivision_operator_norm}
    \|\mathcal{S}_{\boldsymbol{\alpha}^{(\ell)}}\|_\infty=\max\bigg\{\sum_{k\in\mathbb{Z}}|\alpha^{(\ell)}_{2k}|,\; \sum_{k\in\mathbb{Z}}|\alpha^{(\ell)}_{2k+1}|\bigg\}.
\end{equation}
The next result provides a verifiable sufficient condition for the uniform bound~\eqref{eqn:uniform_L}.

\begin{corollary} \label{cor:uniform_stability}
Let $\mathcal{S}_{\boldsymbol{\alpha}}$ be a stationary subdivision operator which is power bounded, that is,
\begin{equation*}
    C_{\boldsymbol{\alpha}}:=\sup_{n\geq0}\|\mathcal{S}_{\boldsymbol{\alpha}}^{\,n}\|_\infty<\infty,
\end{equation*}
and suppose that
\begin{equation}\label{eqn:summable_perturbations}
    \sum_{k=1}^{\infty}\big\|\mathcal{S}_{\boldsymbol{\alpha}^{(k)}}-\mathcal{S}_{\boldsymbol{\alpha}}\big\|_\infty<\infty.
\end{equation}
Then condition~\eqref{eqn:uniform_L} holds with
\begin{equation*}
    L\leq C_{\boldsymbol{\alpha}}\exp\bigg(C_{\boldsymbol{\alpha}}\sum_{k=1}^{\infty}\big\|\mathcal{S}_{\boldsymbol{\alpha}^{(k)}}-\mathcal{S}_{\boldsymbol{\alpha}}\big\|_\infty\bigg).
\end{equation*}
\end{corollary}

\begin{proof}
    Denote $\mathcal{E}_j=\mathcal{S}_{\boldsymbol{\alpha}^{(j)}}-\mathcal{S}_{\boldsymbol{\alpha}}$, $j\in\N$. A direct induction on $q$ verifies the discrete variation-of-constants identity
    \begin{equation*}
        \mathcal{P}_{q,p}=\mathcal{S}_{\boldsymbol{\alpha}}^{\,q-p+1}+\sum_{j=p}^{q}\mathcal{S}_{\boldsymbol{\alpha}}^{\,q-j}\,\mathcal{E}_j\,\mathcal{P}_{j-1,p}, \qquad 1\leq p\leq q,
    \end{equation*}
    with $\mathcal{P}_{p-1,p}=\mathcal{I}$. Indeed, for $q=p$ the identity reads $\mathcal{S}_{\boldsymbol{\alpha}^{(p)}}=\mathcal{S}_{\boldsymbol{\alpha}}+\mathcal{E}_p$, and the induction step follows by writing $\mathcal{P}_{q,p}=(\mathcal{S}_{\boldsymbol{\alpha}}+\mathcal{E}_q)\mathcal{P}_{q-1,p}$. Writing $u_q=\|\mathcal{P}_{q,p}\|_\infty$ and using $\|\mathcal{S}_{\boldsymbol{\alpha}}^{\,n}\|_\infty\leq C_{\boldsymbol{\alpha}}$ for all $n\geq0$, we obtain
    \begin{equation*}
        u_q\leq C_{\boldsymbol{\alpha}}+C_{\boldsymbol{\alpha}}\sum_{j=p}^{q}\|\mathcal{E}_j\|_\infty\, u_{j-1}, \qquad u_{p-1}=1.
    \end{equation*}
    Since $C_{\boldsymbol{\alpha}}\geq1$, the discrete Gr\"onwall inequality yields
    \begin{equation*}
        u_q\leq C_{\boldsymbol{\alpha}}\prod_{j=p}^{q}\big(1+C_{\boldsymbol{\alpha}}\|\mathcal{E}_j\|_\infty\big)\leq C_{\boldsymbol{\alpha}}\exp\bigg(C_{\boldsymbol{\alpha}}\sum_{j=1}^{\infty}\|\mathcal{E}_j\|_\infty\bigg).
    \end{equation*}
    As the bound is independent of $p$, $q$, and $J$, the claim follows.
\end{proof}

Condition~\eqref{eqn:summable_perturbations} is precisely the asymptotic equivalence, in the sense of~\cite{dyn1995analysis}, of the nonstationary family $\{\mathcal{S}_{\boldsymbol{\alpha}^{(k)}} \mid k\in\N\}$ to the stationary scheme $\mathcal{S}_{\boldsymbol{\alpha}}$. It can be verified directly from the masks, since~\eqref{eqn:subdivision_operator_norm}, applied to the difference of the masks, gives
\begin{equation*}
    \big\|\mathcal{S}_{\boldsymbol{\alpha}^{(k)}}-\mathcal{S}_{\boldsymbol{\alpha}}\big\|_\infty=\max\bigg\{\sum_{i\in\mathbb{Z}}\big|\alpha^{(k)}_{2i}-\alpha_{2i}\big|,\; \sum_{i\in\mathbb{Z}}\big|\alpha^{(k)}_{2i+1}-\alpha_{2i+1}\big|\bigg\}.
\end{equation*}
For instance, the masks of the schemes of Examples~\ref{NS example 4-pt} and~\ref{NS example cubic B-spline} converge at a geometric rate to the masks of the $4$-point scheme and of the cubic B-spline scheme, respectively, so that~\eqref{eqn:summable_perturbations} holds; both limit schemes are convergent and hence power bounded, see~\cite{dyn2002subdivision}. Corollary~\ref{cor:uniform_stability} therefore applies to both families, and the corresponding reconstructions are stable uniformly in the number of levels.

Theorem~\ref{thm. stability of reconstruction} ensures that small perturbations in the pyramid representation of a sequence lead to proportionate changes in the subsequent synthesis. This stability property is essential for various applications, including denoising and contrast enhancement, where one typically modifies the detail coefficients before reconstruction. See, for example,~\cite{mattar2023pyramid, mattar2023pseudo}. Conversely, the following theorem shows the stability of the decomposition.

First, we evaluate the operator norm of a decimation operator~\eqref{decimation operator}. Given a decimation operator $\mathcal{D}_{\boldsymbol{\gamma}^{(\ell)}}$ associated with a sequence $\boldsymbol{\gamma}^{(\ell)}$, its operator norm is given by
\begin{equation}\label{eqn:decimation_operator_norm}
    \|\mathcal{D}_{\boldsymbol{\gamma}^{(\ell)}}\|_\infty=\sup_{\|\boldsymbol{c}\|_\infty\leq1}\|\boldsymbol{\gamma}^{(\ell)}*(\boldsymbol{c}\downarrow 2)\|_\infty=\|\boldsymbol{\gamma}^{(\ell)}\|_1.
\end{equation}
Because the elements of $\boldsymbol{\gamma}^{(\ell)}$ sum to one, we have $\|\mathcal{D}_{\boldsymbol{\gamma}^{(\ell)}}\|_\infty\geq 1$ for any such sequence $\boldsymbol{\gamma}^{(\ell)}$. In analogy with~\eqref{eqn:transition_operators}, we introduce the decimation transition operators
\begin{equation}\label{eqn:decimation_transition}
    \mathcal{Q}_{\ell,J}=\mathcal{D}_{\boldsymbol{\gamma}^{(\ell+1)}}\mathcal{D}_{\boldsymbol{\gamma}^{(\ell+2)}}\cdots\mathcal{D}_{\boldsymbol{\gamma}^{(J)}}, \quad \ell=0,\dots,J-1, \qquad \mathcal{Q}_{J,J}=\mathcal{I}.
\end{equation}

\begin{theorem}[Stability of the decomposition]
\label{thm. stability of decomposition}
Let $\boldsymbol{c}^{(J)}$ and $\widetilde{\boldsymbol{c}}^{(J)}$ be two bounded sequences, and let $\{\mathcal{S}_{\boldsymbol{\alpha}^{(\ell)}} \mid \ell=1,\dots,J\}$ be a family of nonstationary subdivision schemes. Denote by $\{\mathcal{D}_{\boldsymbol{\gamma}^{(\ell)}} \mid \ell=1,\dots,J\}$ its corresponding family of decimation operators, associated with sequences $\boldsymbol{\gamma}^{(\ell)}$ whose elements sum to $1$, and let $\mathcal{Q}_{\ell,J}$ be the transition operators~\eqref{eqn:decimation_transition}. Then
\begin{equation}\label{eqn:decomposition_stability_1}
    \|\boldsymbol{c}^{(0)}-\widetilde{\boldsymbol{c}}^{(0)}\|_\infty \leq \|\mathcal{Q}_{0,J}\|_\infty \, \|\boldsymbol{c}^{(J)}-\widetilde{\boldsymbol{c}}^{(J)}\|_\infty,
\end{equation}
and
\begin{equation}\label{eqn:decomposition_stability_2}
    \|\boldsymbol{d}^{(\ell)}-\widetilde{\boldsymbol{d}}^{(\ell)}\|_\infty \leq \|\mathcal{I}-\mathcal{S}_{\boldsymbol{\alpha}^{(\ell)}}\mathcal{D}_{\boldsymbol{\gamma}^{(\ell)}}\|_\infty \, \|\mathcal{Q}_{\ell,J}\|_\infty \, \|\boldsymbol{c}^{(J)}-\widetilde{\boldsymbol{c}}^{(J)}\|_\infty, \quad \ell=1,\dots, J,
\end{equation}
where $\{\boldsymbol{c}^{(0)},\boldsymbol{d}^{(1)},\dots,\boldsymbol{d}^{(J)}\}$ and $\{\widetilde{\boldsymbol{c}}^{(0)},\widetilde{\boldsymbol{d}}^{(1)},\dots,\widetilde{\boldsymbol{d}}^{(J)}\}$ are the multiscale representations~\eqref{NS multiscale transform} of $\boldsymbol{c}^{(J)}$ and $\widetilde{\boldsymbol{c}}^{(J)}$, respectively.
\end{theorem}

\begin{proof}
    Iterating the decimation step of~\eqref{NS multiscale transform} yields $\boldsymbol{c}^{(\ell)}=\mathcal{Q}_{\ell,J}\boldsymbol{c}^{(J)}$ for $\ell=0,\dots,J$, and hence, by the linearity of the operators involved,
    \begin{equation}\label{eqn:c_iterations}
        \|\boldsymbol{c}^{(\ell)}-\widetilde{\boldsymbol{c}}^{(\ell)}\|_\infty = \|\mathcal{Q}_{\ell,J}(\boldsymbol{c}^{(J)}-\widetilde{\boldsymbol{c}}^{(J)})\|_\infty\leq \|\mathcal{Q}_{\ell,J}\|_\infty\|\boldsymbol{c}^{(J)}-\widetilde{\boldsymbol{c}}^{(J)}\|_\infty, \quad \ell=0,\dots,J.
    \end{equation}
    The case $\ell=0$ is~\eqref{eqn:decomposition_stability_1}. Moreover, the detail coefficients satisfy
    \begin{equation*}
    \begin{split}
        \|\boldsymbol{d}^{(\ell)}-\widetilde{\boldsymbol{d}}^{(\ell)}\|_\infty & = \|(\mathcal{I}-\mathcal{S}_{\boldsymbol{\alpha}^{(\ell)}}\mathcal{D}_{\boldsymbol{\gamma}^{(\ell)}})(\boldsymbol{c}^{(\ell)}-\widetilde{\boldsymbol{c}}^{(\ell)})\|_\infty\\
        & \leq \|\mathcal{I}-\mathcal{S}_{\boldsymbol{\alpha}^{(\ell)}}\mathcal{D}_{\boldsymbol{\gamma}^{(\ell)}}\|_\infty \|\boldsymbol{c}^{(\ell)}-\widetilde{\boldsymbol{c}}^{(\ell)}\|_\infty.
    \end{split}
    \end{equation*}
    Therefore, by combining the latter evaluation with~\eqref{eqn:c_iterations} we get~\eqref{eqn:decomposition_stability_2}.
\end{proof}

Naturally, the results of Theorem~\ref{thm. stability of decomposition} agree with Theorem 11 in~\cite{dyn2021linear} when the family of subdivision schemes is stationary. By~\eqref{eqn:decimation_operator_norm} and the submultiplicativity of the operator norm, $\|\mathcal{Q}_{\ell,J}\|_\infty\leq\prod_{m=\ell+1}^{J}\|\boldsymbol{\gamma}^{(m)}\|_1$, so that~\eqref{eqn:decomposition_stability_1} and~\eqref{eqn:decomposition_stability_2} imply the cruder bounds in which the transition norms are replaced by products of the individual norms. The transition-product formulation is sharper. Define
$$ Q_*=\sup_{J\ge1}\max_{0\le\ell\le J}\|Q_{\ell,J}\|_\infty, \qquad R_*=\sup_{\ell\ge1} \|I-S_{\alpha^{(\ell)}}D_{\gamma^{(\ell)}}\|_\infty. $$
If \(Q_*<\infty\) and \(R_*<\infty\), then \eqref{eqn:decomposition_stability_1} and \eqref{eqn:decomposition_stability_2} give decomposition stability uniformly in both \(J\) and \(\ell\). These conditions may hold even when the corresponding products of individual operator norms grow. In practice, the bound~\eqref{eqn:decomposition_stability_1} is also affected by the truncation threshold $\varepsilon$: a faster decay of the elements of $\boldsymbol{\gamma}^{(\ell)}$ permits, for a given $\varepsilon$, a truncation with a smaller error. We note, however, that a faster decay does not by itself imply a smaller $\ell^1$-norm. In the particular case where the multiscaling~\eqref{NS multiscale transform} is interpolating, $\boldsymbol{\gamma}^{(\ell)}=\boldsymbol{\delta}$ and $\|\mathcal{D}_{\boldsymbol{\gamma}^{(\ell)}}\|_\infty=1$. This is the lowest value possible for this operator; in this case, the operators $\mathcal{Q}_{\ell,J}$ reduce to plain downsampling with $\|\mathcal{Q}_{\ell,J}\|_\infty=1$, which may lead to equality in~\eqref{eqn:decomposition_stability_1}.

Another valuable insight stemming from~\eqref{eqn:decomposition_stability_2} is the following. First, observe the relation
\begin{equation*}
    \|[\mathcal{I}-\mathcal{S}_{\boldsymbol{\alpha}^{(\ell)}}\mathcal{D}_{\boldsymbol{\gamma}^{(\ell)}}(\cdot)]\downarrow 2\|_\infty \leq \|\mathcal{I}-\mathcal{S}_{\boldsymbol{\alpha}^{(\ell)}}\mathcal{D}_{\boldsymbol{\gamma}^{(\ell)}}\|_\infty,
\end{equation*}
where the operator on the left hand side returns the even detail coefficients when applied to a sequence, and its operator norm equals the residual of the convolutional equation~\eqref{convolutional equation for finding gamma}, namely $\|\boldsymbol{\delta}-(\boldsymbol{\alpha}^{(\ell)}\downarrow2)*\boldsymbol{\gamma}^{(\ell)}\|_1$, which measures how well $\boldsymbol{\gamma}^{(\ell)}$ and $\boldsymbol{\alpha}^{(\ell)}$ satisfy the equation. This identity follows from Theorem 4 in~\cite{dyn2021linear} and its proof. To quantify the effect of the truncation, denote by $\boldsymbol{\gamma}^{(\ell)}_\varepsilon$ the sequence obtained from $\boldsymbol{\gamma}^{(\ell)}$ by discarding the elements smaller in magnitude than the threshold $\varepsilon$. Then, the residual associated with the truncated sequence satisfies
\begin{equation*}
    \|\boldsymbol{\delta}-(\boldsymbol{\alpha}^{(\ell)}\downarrow2)*\boldsymbol{\gamma}^{(\ell)}_\varepsilon\|_1 = \|(\boldsymbol{\alpha}^{(\ell)}\downarrow2)*(\boldsymbol{\gamma}^{(\ell)}-\boldsymbol{\gamma}^{(\ell)}_\varepsilon)\|_1 \leq \|\boldsymbol{\alpha}^{(\ell)}\downarrow2\|_1\,\|\boldsymbol{\gamma}^{(\ell)}-\boldsymbol{\gamma}^{(\ell)}_\varepsilon\|_1,
\end{equation*}
where $\|\boldsymbol{\gamma}^{(\ell)}-\boldsymbol{\gamma}^{(\ell)}_\varepsilon\|_1$, the sum of the magnitudes of the discarded elements, is nondecreasing in $\varepsilon$; the latter holds provided that the retained coefficient sum is bounded away from zero. Hence, larger truncation thresholds may result in even-indexed detail coefficients that deviate further from zero, an effect that is controlled linearly by the truncation error.

\section{Geometry-driven pyramids and their applications} \label{sec:geom}

Until now, we have demonstrated how to construct the pyramid using a nonstationary upsampling operator, even in the absence of guaranteed interpolation. However, the practical relevance of this construction may not yet be fully apparent. In this section, we present several concrete constructions that highlight the applicability of the framework from a geometric perspective.

All figures and examples presented in this section are available as open-source code to support reproducibility and can be accessed at: \url{https://github.com/HadarLandau/NSP}.

\subsection{Reproducing exponential polynomials: The refinement}

We borrow similar notions from~\cite{conti2011algebraic} and start by defining spaces of exponential polynomials. Let $T \in \mathbb{Z}_+$, and let $\boldsymbol{\lambda}=(\lambda_0,\lambda_1, \dots, \lambda_T)$ be a vector of real coefficients with $\lambda_T \neq 0$. Let $D^n=d^n/dx^n$ denote the $n$-th order
differentiation operator. Let $V_{T,\boldsymbol{\lambda}}$ be the space of exponential polynomials given by
\begin{equation}
    V_{T,\boldsymbol{\lambda}} = \bigg\{ f \colon \mathbb{R} \rightarrow \mathbb{C} \, \big| \, \sum_{j=0}^T \lambda_j D^j f=0 \bigg\},
\end{equation}
Being the solution space of a linear differential equation with constant coefficients, $V_{T,\boldsymbol{\lambda}}$ is a $T$-dimensional subspace of $C^\infty(\mathbb{R},\mathbb{C})$.

Now, define the function $\boldsymbol{\lambda}(z)=\sum^T_{j=0} \lambda_j z^j$ where $z\in\mathbb{C}$ (often called the $z$-transform of $\boldsymbol{\lambda}$), and let $\{(\zeta_l,\tau_l)\}_{l=1,\dots,N}$ denote the zeros of $\boldsymbol{\lambda}(z)$ together with their multiplicities. Namely,
\[
     \boldsymbol{\lambda}^{(r)}(\zeta_l)=0, 
\]
where $r=0,\dots,\tau_l -1$, $l=1,\dots,N$, $\tau_l \in \mathbb{N}$, and $T=\textstyle\sum_{l=1}^N \tau_l$. Under these notations, the subspace $V_{T,\boldsymbol{\lambda}}$ reduces to
\begin{equation}
    V_{T,\boldsymbol{\lambda}}=\operatorname{span} \{x^r e^{\zeta_l x} \;|\; r=0,\dots,\tau_l -1, \ l=1,\dots,N \}.
\end{equation}

In the following, we provide the definition of a subdivision scheme that reproduces the space of exponential polynomials. 

\begin{definition} [Exponential reproduction] \label{def:exponential_reproduction}
Let $\boldsymbol{\mathcal{S}}=\{\mathcal{S}_{\boldsymbol{\alpha}^{(k)}} \mid k\in \N \}$ be a convergent scheme. We say that $\boldsymbol{\mathcal{S}}$ is \emph{exponential reproducing} with respect to $V_{T,\boldsymbol{\lambda}}$ if, for any function $f \in V_{T,\boldsymbol{\lambda}}$, any starting level $m \in \N$, and the initial sequence $\boldsymbol{c}^{(0)}$ whose elements are the samples $c^{(0)}_j=f\big(2^{-(m-1)}j\big)$, $j\in\mathbb{Z}$, it holds that
\begin{equation}
    \lim_{k\rightarrow\infty} \, \sup_{j\in\mathbb{Z}} \, \Big| f\big(2^{-(m+k)}j\big) - \big( \mathcal{S}_{\boldsymbol{\alpha}^{(m+k)}} \mathcal{S}_{\boldsymbol{\alpha}^{(m+k-1)}} \cdots \mathcal{S}_{\boldsymbol{\alpha}^{(m)}} \boldsymbol{c}^{(0)} \big)_j \Big| = 0.
\end{equation}
\end{definition}
Note that reproduction is required for every starting level $m$, since the operators of a nonstationary scheme vary with the level. For noninterpolatory schemes, reproduction is understood with respect to a suitable parametrization (shift) of the refined grids; see~\cite{conti2011algebraic} for the general, shift-dependent formulation.

We recall an example from \cite{conti2011algebraic} of a nonstationary subdivision scheme that reproduces the space
\begin{equation*}
    V_{T,\boldsymbol{\lambda}} = \text{span}\{1, x, e^{tx}, e^{-tx}\},
\end{equation*}
which is suitable for representing conic sections. This space is obtained by taking $T=4$ and $\boldsymbol{\lambda}=(0, 0, -t^2,0, 1)$ for some $t\in\mathbb{R}$, resulting in $\boldsymbol{\lambda}(z)=z^4-t^2z^2$. Depending on the value of the parameter $t$, this space captures different types of conics. For example, hyperbolic functions are obtained for $t=1$ and trigonometric functions for $t=i$, while for $t=0$ the zero $z=0$ has multiplicity four and the space degenerates to the space of cubic algebraic polynomials.

\begin{example} \label{example:geometric_scheme}
The scheme $\boldsymbol{\mathcal{S}}=\{\mathcal{S}_{\boldsymbol{\alpha}^{(k)}} \mid k\in\N\}$ is given by the refinement rules, written here as mapping level $k$ to level $k+1$, which matches the convention of~\eqref{nonstationary refinement rule} up to an index shift:
\begin{equation} \label{NS geometric scheme}
\begin{split}    
    & c_{2j}^{(k+1)}=\frac{a^{(k)}}{4(v^{(k)}+1)} (c_{j-2}^{(k)}+c_{j+2}^{(k)})+\frac{1+2v^{(k)}(b^{(k)}+2a^{(k)})}{4(v^{(k)}+1)}(c_{j-1}^{(k)}+c_{j+1}^{(k)}) \\
    & \qquad \qquad +\frac{4v^{(k)}(1-b^{(k)}-2a^{(k)})-2a^{(k)}+2}{4(v^{(k)}+1)}c_{j}^{(k)},\\
    & c_{2j+1}^{(k+1)}=\frac{2a^{(k)}(v^{(k)}+1)+b^{(k)}}{4(v^{(k)}+1)}(c_{j-1}^{(k)}+c_{j+2}^{(k)})\\
    & \qquad \qquad +\frac{(2-2a^{(k)})(v^{(k)}+1)-b^{(k)}}{4(v^{(k)}+1)}(c_{j}^{(k)}+c_{j+1}^{(k)}),
\end{split}
\end{equation}
where 
\begin{equation}
\begin{split}
    a^{(k)} & =\frac{(2+\sqrt{2(v^{(k)}+1)})(2-v^{(k)}\sqrt{2(v^{(k)}+1)})}{8v^{(k)}(v^{(k)}-1)\sqrt{2(v^{(k)}+1)}(v^{(k)}+3+2\sqrt{2(v^{(k)}+1)})},\\
    b^{(k)} & =\frac{(v^{(k)}+1)(v^{(k)}-2)-2\sqrt{2(v^{(k)}+1)}}{2v^{(k)}\sqrt{2(v^{(k)}+1)}(v^{(k)}+3+2\sqrt{2(v^{(k)}+1)})},
\end{split}
\end{equation}
with $v^{(k)}=\sqrt{\frac{1+v^{(k-1)}}{2}}, \ k \geq 0$ and $v^{(-1)}>-1$. For $v^{(k)}=1$, which occurs in the purely polynomial case below, the expressions for $a^{(k)}$ and $b^{(k)}$ are of the indeterminate form $0/0$, and are understood as their limit values $a^{(k)}=-5/64$ and $b^{(k)}=-3/16$, obtained by continuity. We also note that, since $v^{(k)}\rightarrow1$ as $k$ grows, the direct evaluation of $a^{(k)}$ is prone to numerical cancellation at high refinement levels; this is immaterial in our experiments, where only a few levels are applied.

The determination of the initial value $v^{(-1)}$, is carried out under the assumption that the initial data are obtained by uniformly sampling a function $f \in V_{T,\lambda}$. The value of $v^{(-1)}$ is then chosen as
\begin{equation}
    v^{(-1)}=
    \begin{cases}
        1 & \text{when $f$ is purely polynomial,}\\
        \cosh(\sigma) & \text{when $f$ is hyperbolic,}\\
        \cos(\sigma) & \text{when $f$ is trigonometric,}
    \end{cases}
\end{equation}
where $\sigma$ is a positive number that denotes the spacing between the uniformly sampled points. For instance, if the unit circle is sampled at $N$ equispaced points, then $\sigma=2\pi/N$.
\end{example}

\begin{figure}[htbp]
    \centering  
    \begin{subfigure}{0.35\textwidth}
         \centering         \includegraphics[width=\textwidth]{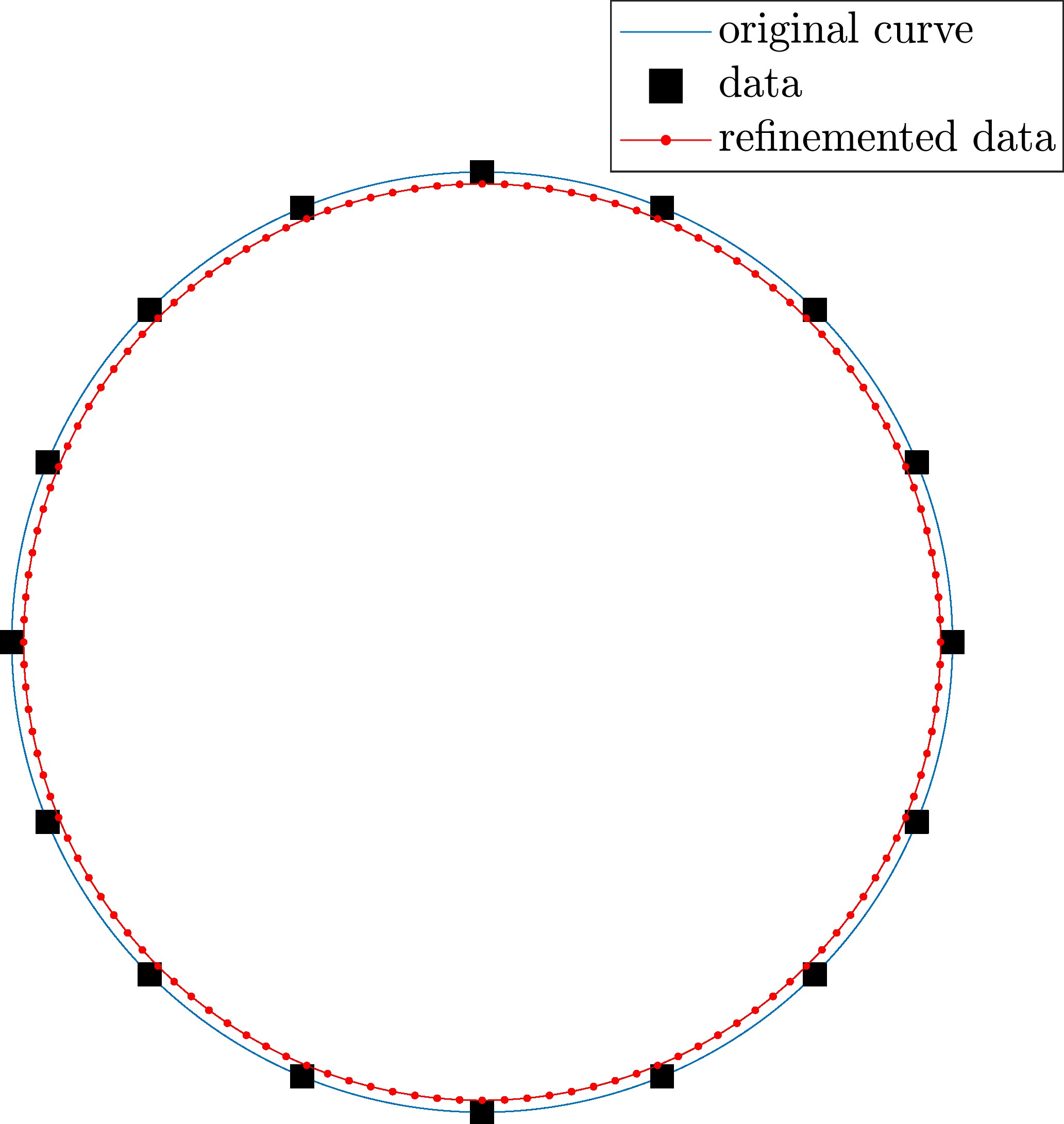}
    \end{subfigure}
    \qquad \qquad
    \begin{subfigure}{0.35\textwidth}
         \centering         \includegraphics[width=\textwidth]{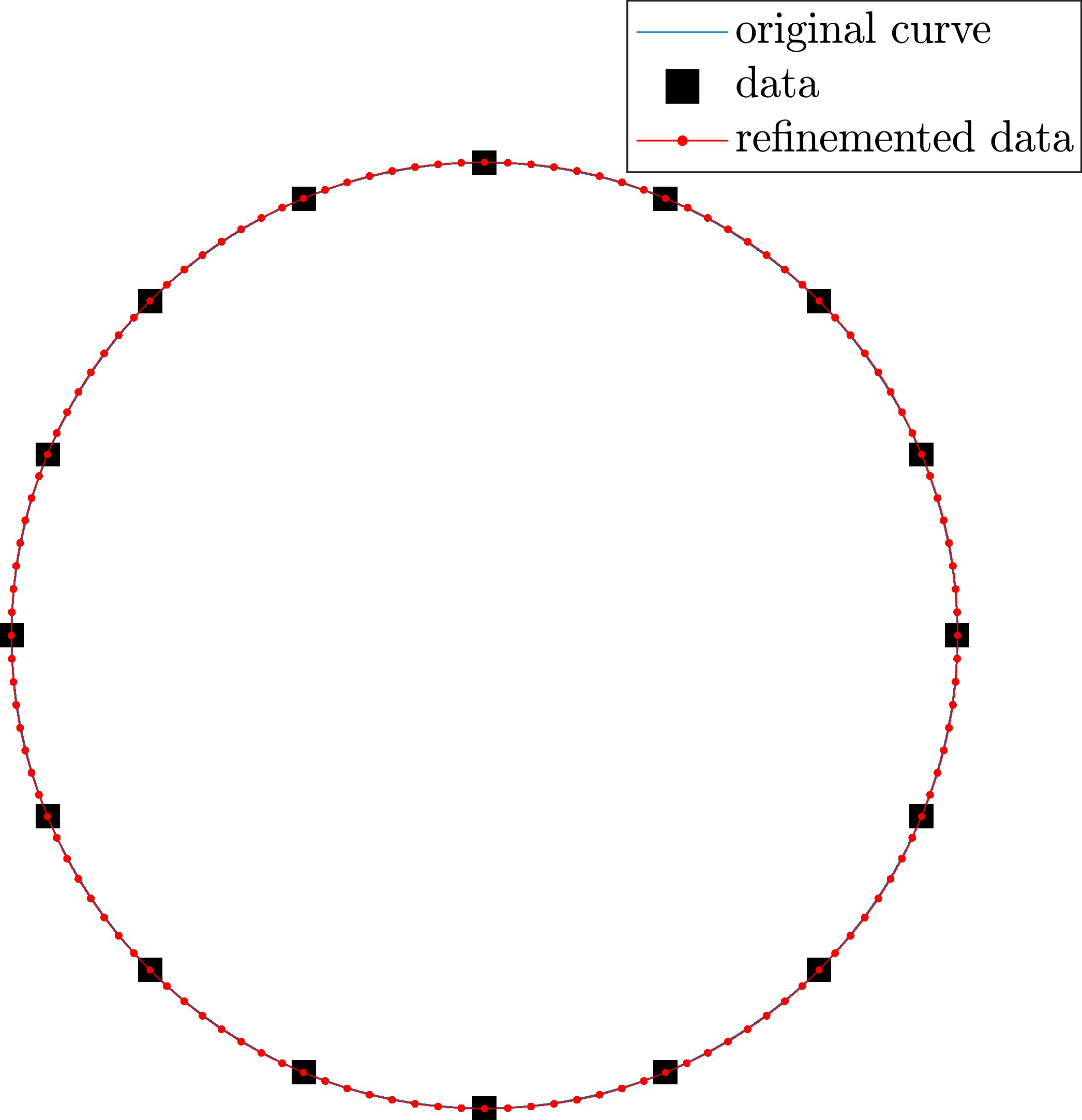}
     \end{subfigure}
    \caption{A comparison is shown between the refinement results for a circular shape obtained using two subdivision schemes. The original curves (blue), the initial coarse points (black squares), and the refined points (red dots) after three refinement steps are displayed. The left plot, using the cubic B-spline scheme, yields an inaccurate approximation of the circle, whereas the right plot, using the geometric scheme \eqref{NS geometric scheme}, matches the original curve up to the machine precision, thereby highlighting its exponential reproduction property.}
    \label{fig:reproducing circles}
\end{figure}

The exponential reproduction property of the scheme of Example~\ref{example:geometric_scheme} is demonstrated in Figure~\ref{fig:reproducing circles}. Note that for this example, only the finite-depth estimate is being invoked. We begin by sampling the unit circle, an example of a trigonometric polynomial, at nine equally spaced points. These samples are then refined over three levels using two schemes: the cubic B-spline scheme and the scheme of Example~\ref{example:geometric_scheme}. As can be observed, the latter successfully reproduces the circle, keeping the refined points precisely on the unit circle, whereas the cubic B-spline scheme yields points that drift away from it.

\subsection{Multiscale analysis of perturbed circles}

Due to the locality of the subdivision scheme~\eqref{NS geometric scheme} and its ability to reproduce conic sections, its associated nonstationary pyramid transform provides a potential tool for evaluating the circularity of a given set of samples. In particular, small detail coefficients indicate that the observed data points are consistent with the space of curves reproduced by the scheme, a model that contains the circles; see Remark~\ref{rem:ellipses} below for the precise scope of this criterion. Within the setting of the following experiments, where the data are perturbations of uniformly sampled circles, the detail coefficients thus quantify how closely the sampled points resemble a circle. In the following series of examples, we illustrate the practical validity of the proposed application. We generate several sets of samples drawn from the unit circle, each corrupted by varying levels and types of oscillatory noise, and analyze their circular structure, as reflected by the first four levels of their multiscale representations. We quantify the overall circularity of a given data sequence $\boldsymbol{c}$ at level $\ell$ by defining $\nu_\ell(\boldsymbol{c})$ as the average norm of its detail coefficients at level $\ell$, given by
\begin{equation} \label{eqn:nu}
    \nu_\ell(\boldsymbol{c}) = \frac{1}{N_{\ell}}\sum_{i=1}^{N_{\ell}}\big\|d_i^{(\ell)}\big\|_2,
\end{equation}
where $N_{\ell}$ is the number of detail coefficients at the $\ell$-th level, and $\|\cdot\|_2$ is the Euclidean norm in $\mathbb{R}^2$.

\begin{remark} \label{rem:implementation}
While the theory of Section~\ref{sec:pyramid} is formulated for scalar bi-infinite sequences, the planar data considered here are processed componentwise; that is, the transform~\eqref{NS multiscale transform} is applied separately to each coordinate, and each detail coefficient is a planar vector. Since the analyzed curves are closed, the sequences are treated as periodic, and all operators are applied with cyclic indexing. In this setting, Theorem~\ref{thm. details decay} applies to the coordinate functions of a smooth $2\pi$-periodic parametrization of the curve.
\end{remark}

\begin{remark}[Scope of the circularity criterion] \label{rem:ellipses}
Exponential reproduction gives a one-way implication. With the exact
decimation operators and the appropriate parametrization, uniform samples
of a curve belonging componentwise to the space reproduced by the scheme
of Example~4.2 yield zero detail coefficients. When the truncated and
renormalized decimation operators are used, the corresponding details are
zero up to truncation and numerical errors. The converse need not hold:
vanishing details over finitely many decomposition levels show that the data
lie in the range of the corresponding finite-depth subdivision operators, but
do not by themselves establish membership in the reproduced
exponential-polynomial space.

With the trigonometric initialization, the reproduced space contains, in
each coordinate, the span
\[
    \operatorname{span}\{1,x,\cos(\sigma x),\sin(\sigma x)\}.
\]
Consequently, uniform samples of any planar curve of the form
\[
    \boldsymbol{c}(t)
    =
    \boldsymbol{p}+\boldsymbol{q}t+
    A
    \begin{pmatrix}
        \cos(\sigma t)\\
        \sin(\sigma t)
    \end{pmatrix},
    \qquad
    \boldsymbol{p},\boldsymbol{q}\in\mathbb{R}^{2},
    \quad A\in\mathbb{R}^{2\times2},
\]
yield zero detail coefficients. Periodicity implies $\boldsymbol{q}=\boldsymbol{0}$. This family therefore includes all ellipses centered at $\boldsymbol{p}$, including degenerate ones, and not only circles.

Accordingly, $\nu_\ell$ in~\eqref{eqn:nu} should be interpreted
as a refinement-residual diagnostic rather than as a certificate of
circularity or a metric distance to the reproduced trigonometric model.
Circles may be distinguished within that model by an additional test, for
example by verifying that the sample radii relative to their centroid are
approximately constant, together with the prescribed uniform ordering of
the samples. In the experiments below, the data are generated as
perturbations of uniformly sampled circles. Within this controlled setting,
$\nu_\ell$ provides a useful relative measure of departure from the
circle-reproducing refinement model, and we therefore refer to it as a
circularity measure. Finally, $\nu_\ell$ is scale-dependent, since
\[
    \nu_\ell(a\boldsymbol{c}+\boldsymbol{p})
    =
    |a|\,\nu_\ell(\boldsymbol{c}).
\]
\end{remark}

We proceed by deriving the pyramid transform that corresponds to the subdivision scheme appearing in Example~\ref{example:geometric_scheme}. To this end, we substitute $\boldsymbol{\alpha}_1^{(\ell)}$ into \eqref{convolutional equation for finding gamma} to determine the decimation operator for each level. Note that, in the subdivision scheme, \( v^{(\ell)} \) is computed to match the number of samples, which doubles each time \( \ell \) increases. In contrast, during the decomposition process, the number of samples is halved at each level. Since the precise choice of \( v^{(\ell)} \) is critical for preserving the reproducing properties of the scheme, this mismatch is handled by selecting, at each level $\ell$, a new initial value $v^{(-1)}$ that corresponds to the number of samples present at that level, and by taking $\boldsymbol{\alpha}_1^{(\ell)}$ to be the mask of the corresponding first refinement step of the scheme of Example~\ref{example:geometric_scheme}. As before, we denote the $\ell$-th level decimation operator that we calculate in~\eqref{convolutional equation for finding gamma} by $\boldsymbol{\gamma}_1^{(\ell)}$. We truncate the decimation with $\varepsilon=10^{-15}$ which gave rise to $33$ nonzero coefficients and normalize each $\boldsymbol{\gamma}_1^{(\ell)}$. Both $\boldsymbol{\alpha}_1^{(\ell)}$ and $\boldsymbol{\gamma}_1^{(\ell)}$ are substituted into \eqref{NS multiscale transform} to compute $\boldsymbol{c}^{(\ell-1)}$ and $\boldsymbol{d}^{(\ell)}$.

\begin{figure}[htbp]
    \centering
    \begin{subfigure}[b]{0.35\textwidth}
         \centering         \includegraphics[width=.9\textwidth]{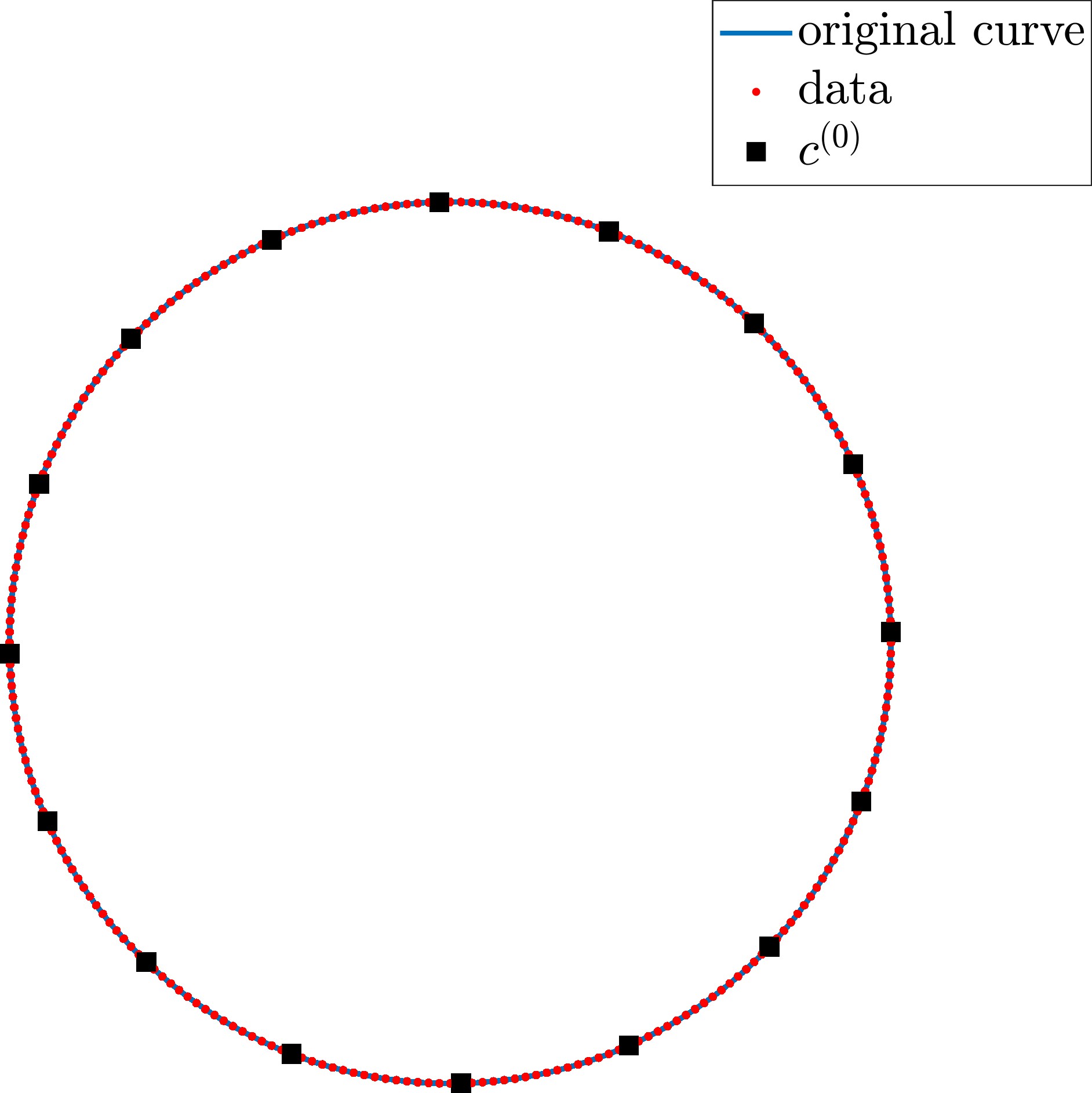}
         \label{fig:a1_circle}
     \end{subfigure} \qquad
     \begin{subfigure}[b]{0.3\textwidth}
         \centering         \includegraphics[width=\textwidth]{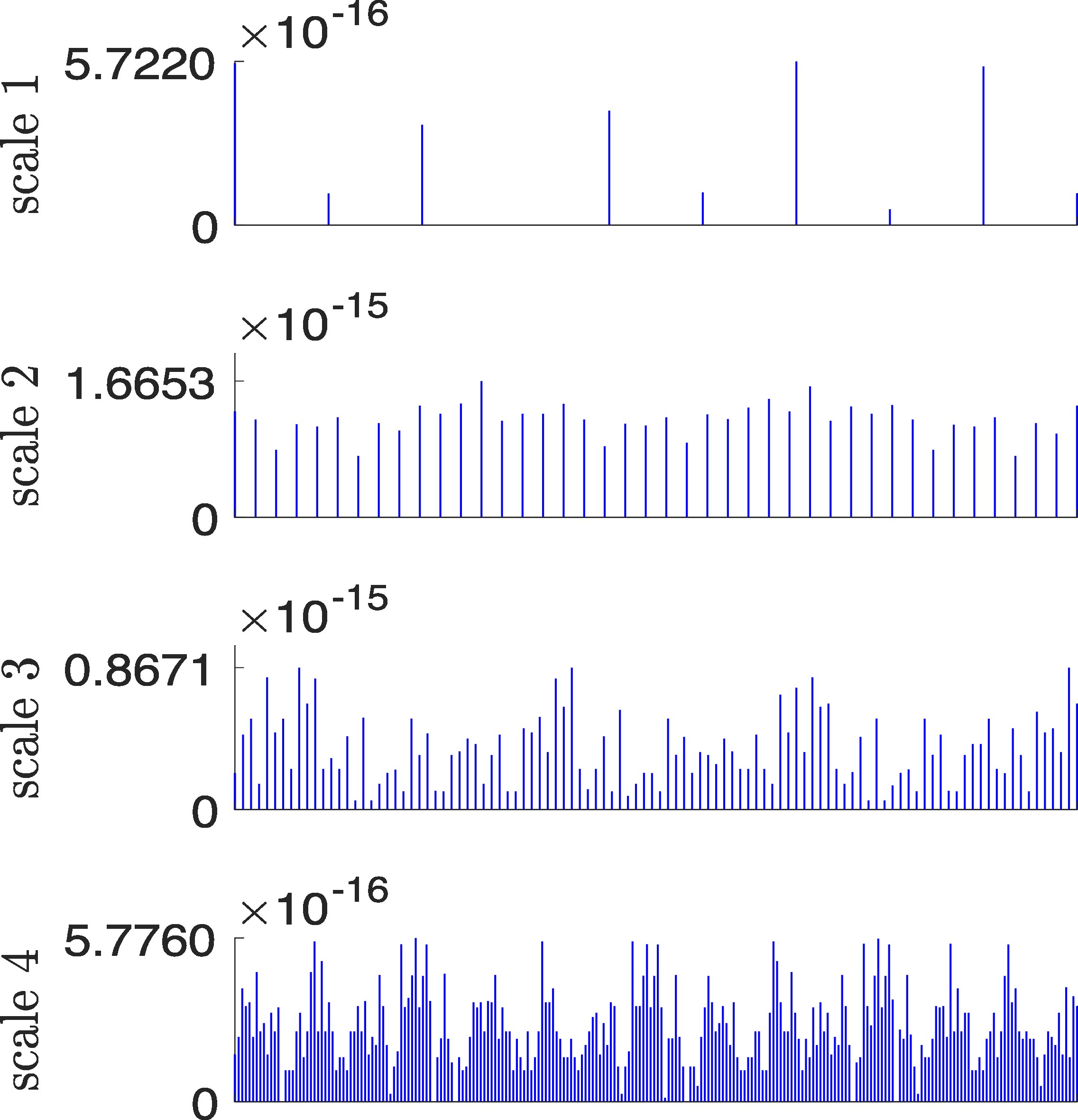}
   \label{fig:a1_circle_dtls}
    \end{subfigure}
     \caption{Points on the unit circle and their multiscale transform. On the left, $256$ equispaced samples of the unit circle (in red) are shown together with $16$ data points (black squares) obtained by applying the analysis \eqref{NS multiscale transform} four times. On the right, the Euclidean norms of the detail coefficients are displayed. Their near-zero values confirm the reproduction of the sampled circle by the chosen scheme.}
     \label{fig:clean circle}
\end{figure}

\begin{figure}[htbp]
    \centering
    \begin{subfigure}{0.32\textwidth}
         \centering         \includegraphics[width=\textwidth]{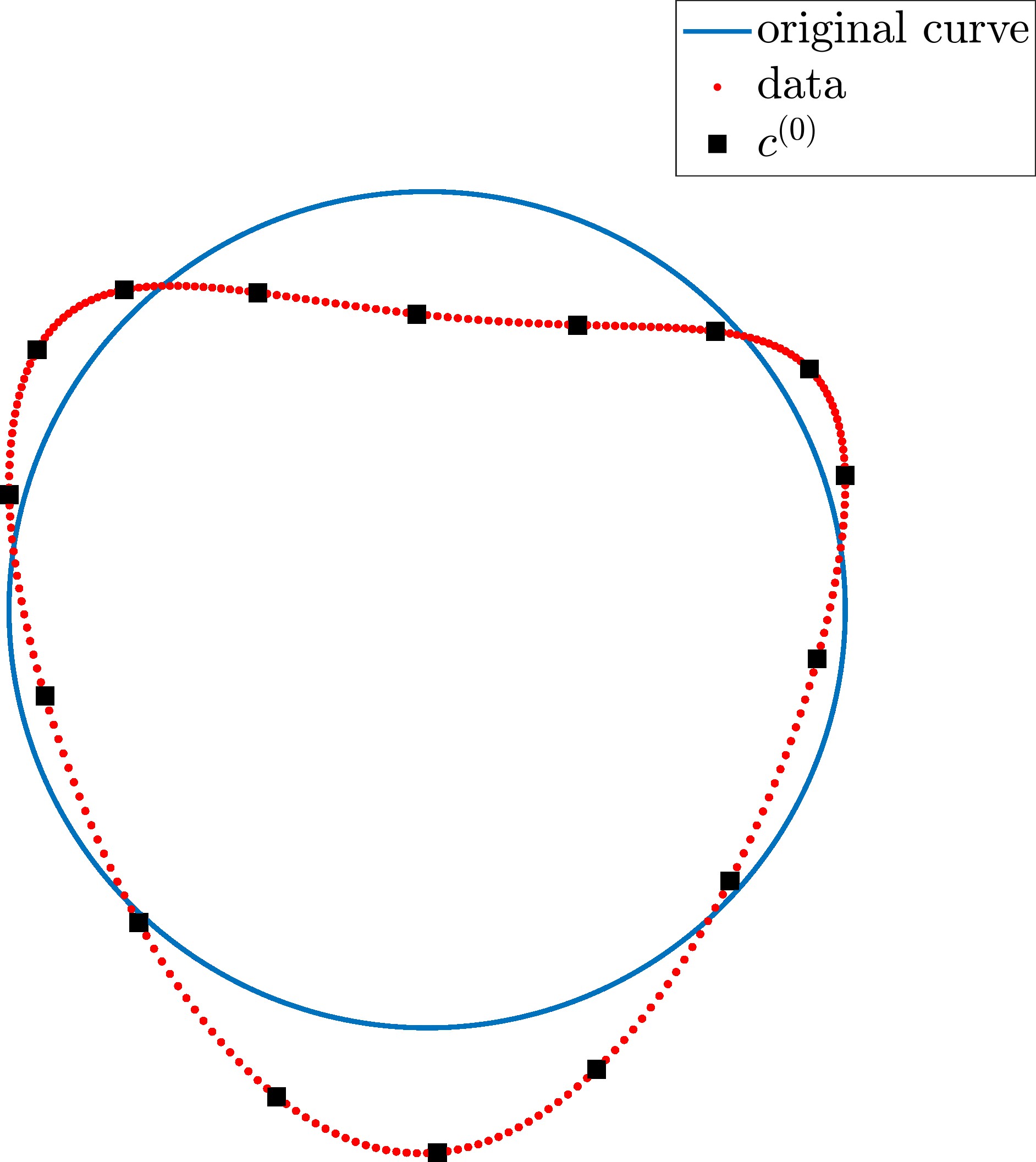}
     \end{subfigure}
     \begin{subfigure}{0.32\textwidth}
         \centering         \includegraphics[width=\textwidth]{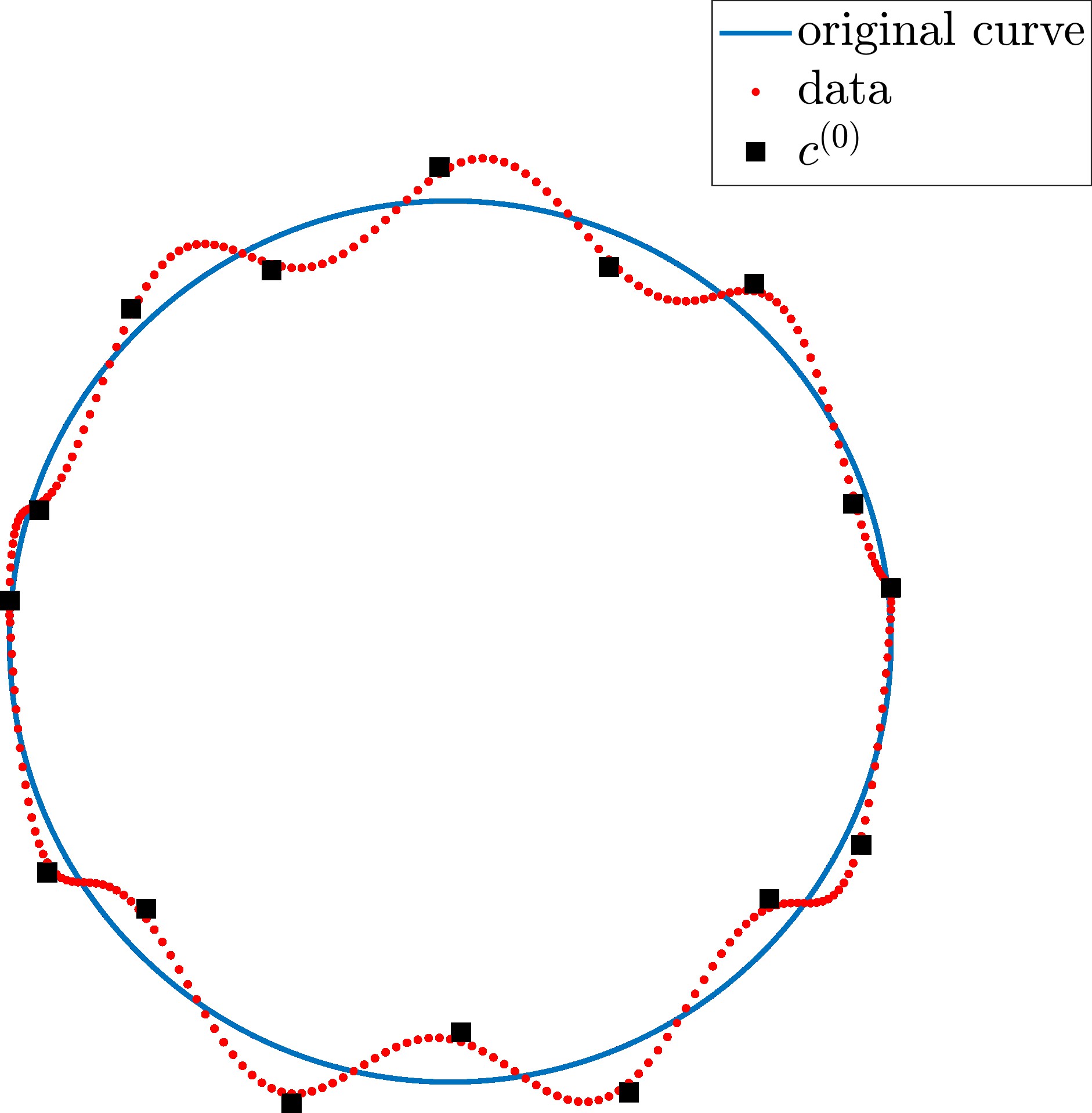}
    \end{subfigure}
    \begin{subfigure}{0.32\textwidth}
         \centering         \includegraphics[width=\textwidth]{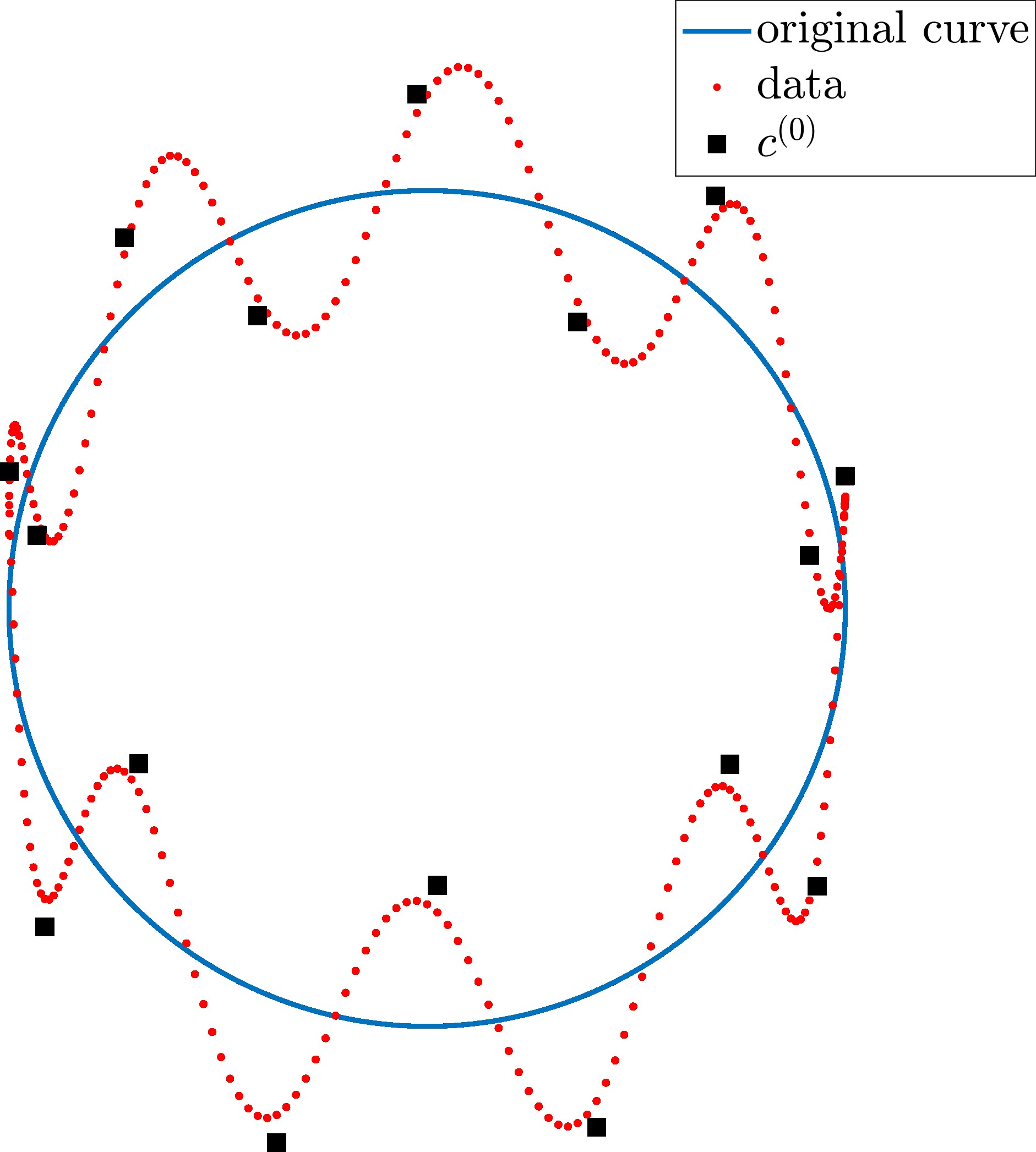}
    \end{subfigure} \\
     \begin{subfigure}{0.32\textwidth}
         \centering        \includegraphics[width=\textwidth]{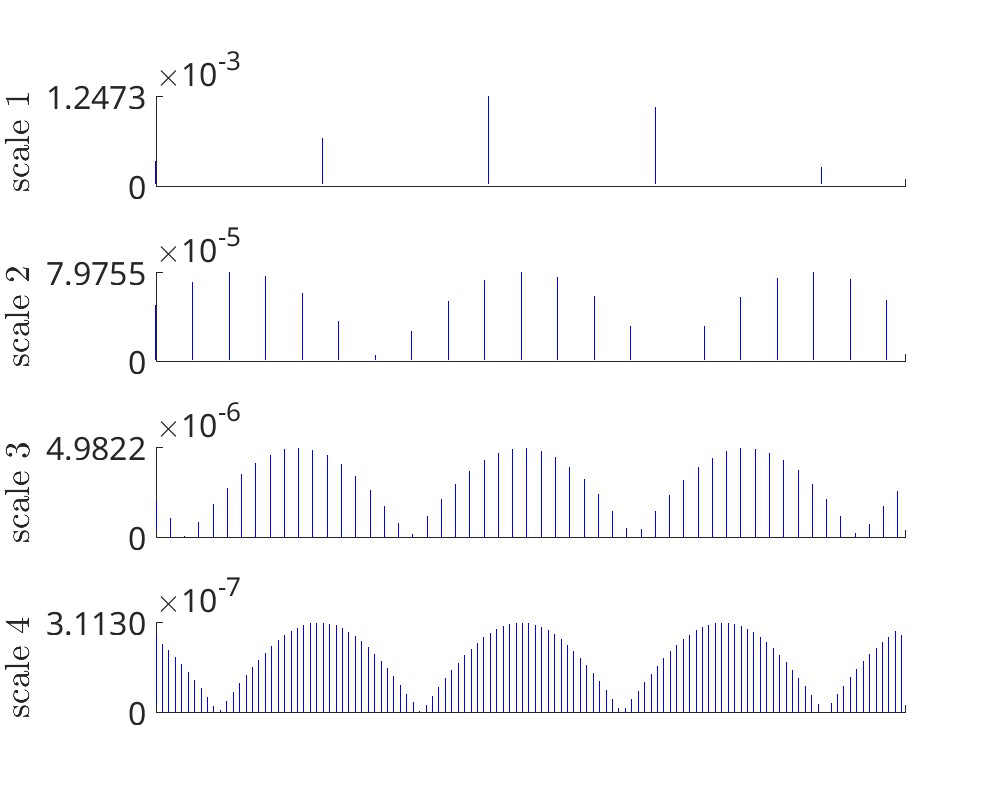}
    \end{subfigure} 
    \begin{subfigure}{0.32\textwidth}
         \centering         \includegraphics[width=\textwidth]{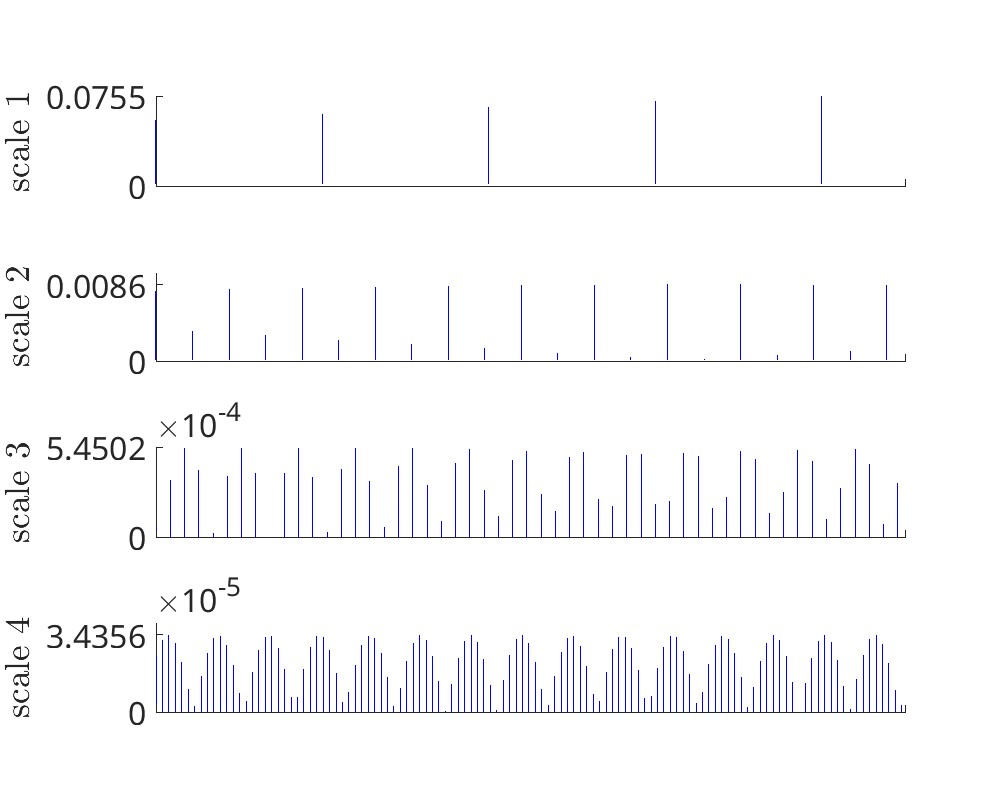}
     \end{subfigure}
    \begin{subfigure}{0.32\textwidth}
         \centering         \includegraphics[width=\textwidth]{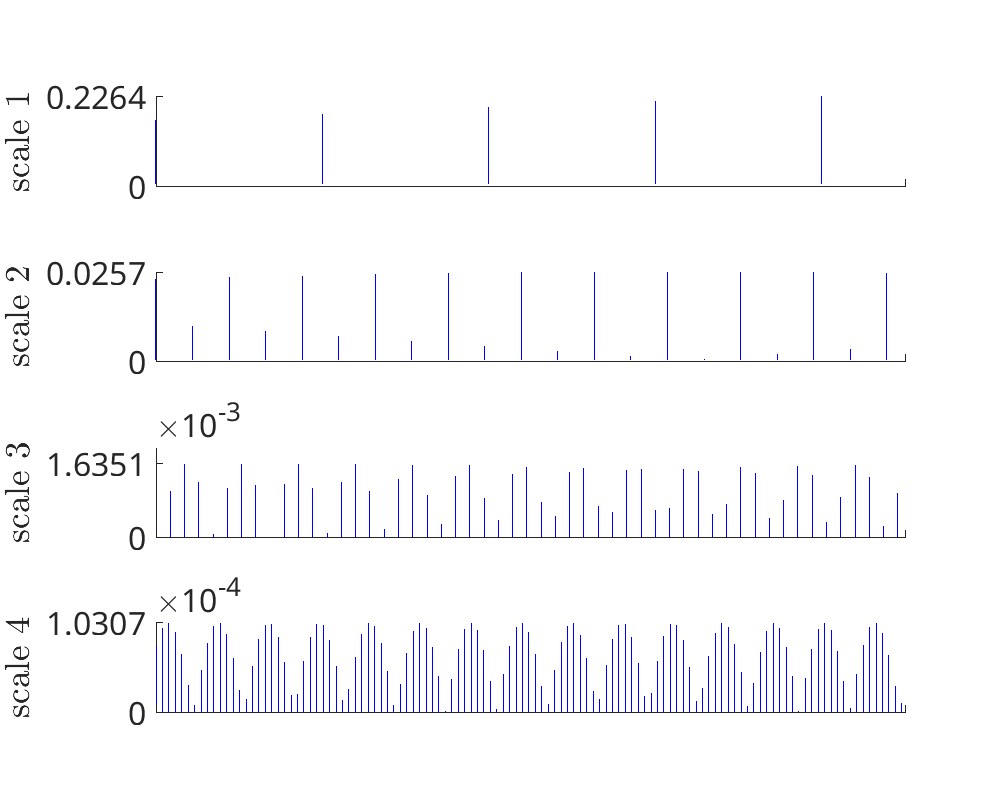}
     \end{subfigure}
     \caption{Quantifying similarity to a circular shape using multiscale transforms. Three examples of points sampled from perturbed circular curves (``wavy circles'') are shown. For each example, the upper plot displays 256 points along the curve and 16 coarse points (black squares) obtained after four decomposition steps of the transform. The lower plots illustrate the Euclidean norms of the corresponding detail coefficients. From left to right, perturbations increase, leading to larger detail coefficients. Smaller detail coefficients indicate greater similarity to a perfect circle.}
     \label{fig:wavy_circle_pyramid}
\end{figure}

\begin{figure}[htbp]
    \centering    
    \includegraphics[width=0.46\linewidth]{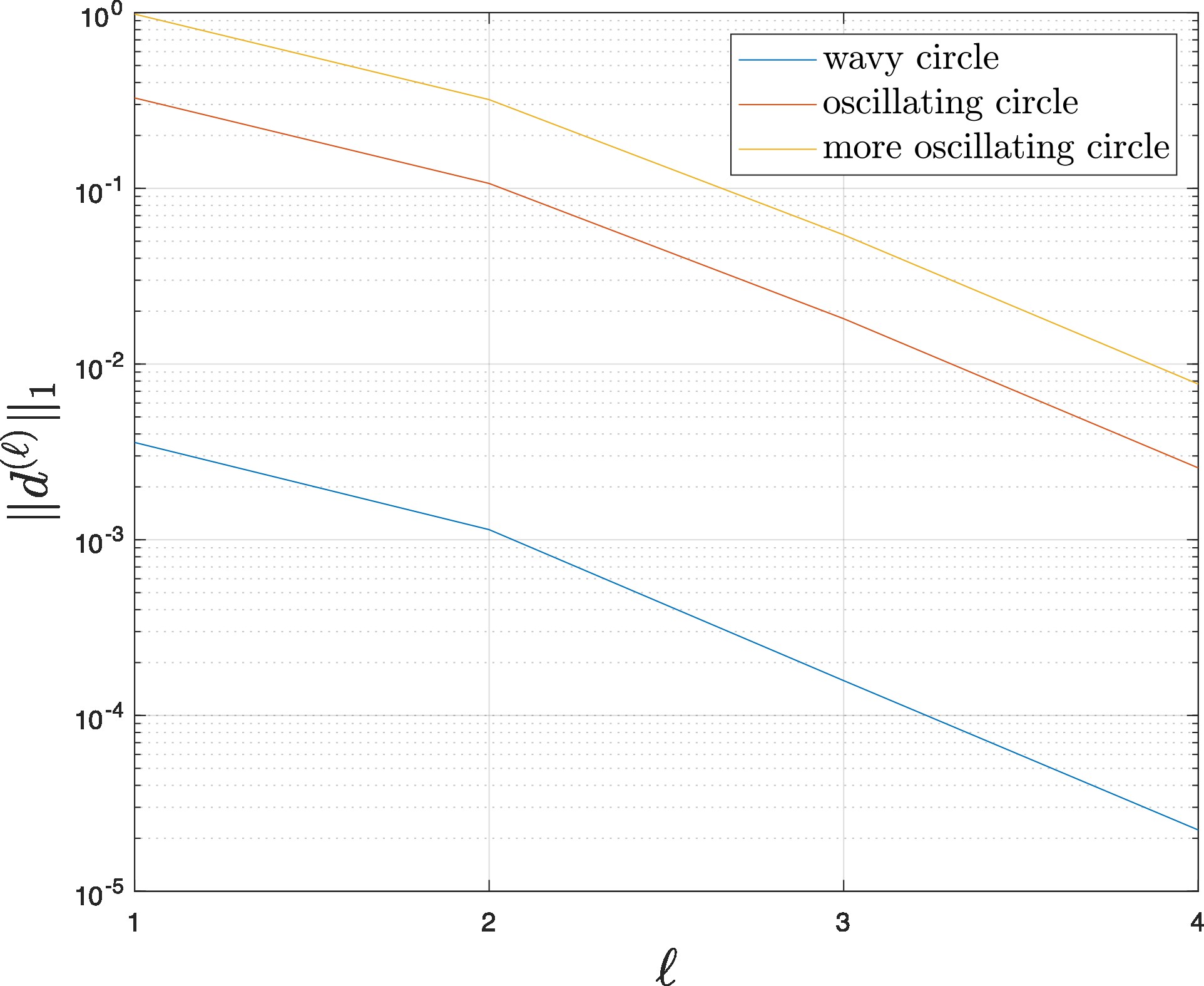} \quad 
    \includegraphics[width=0.42\linewidth]{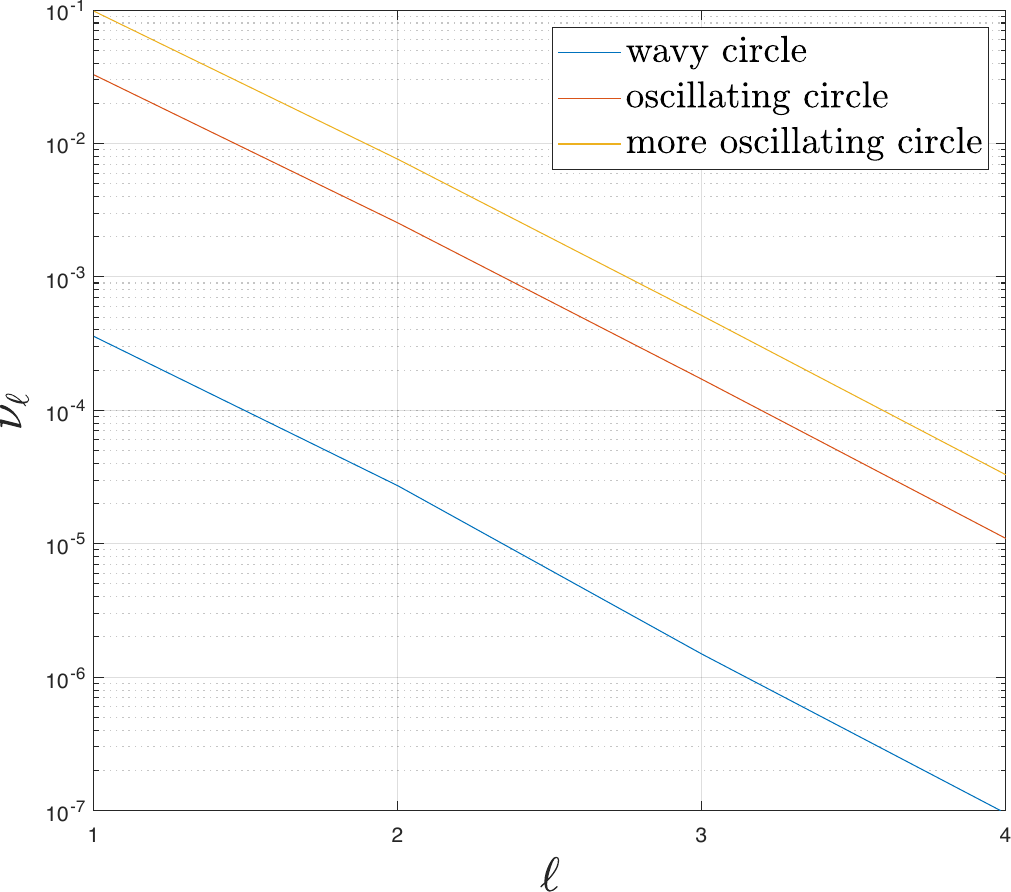}
    \caption{Log-scale plots of detail coefficients decay for three example shapes shown in \figref{fig:wavy_circle_pyramid}. The blue curve (``wavy circle'') corresponds to the leftmost shape, the red curve (``oscillating circle'') to the central shape, and the yellow curve (``more oscillating circle'') to the rightmost shape. The left panel shows the decay of the detail coefficients, measured by their $\ell^1$-norm, indicating their magnitude and sparsity for each level $\ell=1,\dots,4$. The right panel shows the decay, measured by $\nu_\ell(\boldsymbol{c})$, which quantifies the average detail magnitude for each dataset across scales $\ell=1,\dots,4$. Together, these two perspectives illustrate the hierarchical relationships revealed by the pyramid transform across the three examples.
}
    \label{fig:circle_detection_log}
\end{figure}

In our first example shown in \figref{fig:clean circle}, we present a clean planar circle, sampled at $256$ uniformly spaced points, along with its multiscale decomposition. The analysis process yields a coarse representation of $16$ data points, which are also marked in this figure. On the right, the norms of the corresponding detail coefficients are presented, as we recall that each coefficient is a planar vector.
For this uncorrupted data we obtain $\nu_\ell(\boldsymbol{c}_{\text{clean}}) \approx 0$ (up to machine precision) across all levels $\ell=1,\dots,4$, quantitatively confirming the exponential reproduction property of the circle-reproducing scheme. 
Next, in \figref{fig:wavy_circle_pyramid} we illustrate three examples of unit circles perturbed by oscillations of increasing magnitude from left to right. The bottom row shows the corresponding multiscale pyramid analysis for each case. As the magnitude of perturbations grows, the shape deviates further from a perfect circle, resulting in larger values of $\nu_\ell(\boldsymbol{c})$ across all levels. Lastly, in \figref{fig:circle_detection_log} we provide comparison graphs for the $\ell^1$-norm and $\nu_\ell(\boldsymbol{c})$ as a function of the level $\ell$. This comparison clearly demonstrates a consistent increase in both metrics as shape perturbations grow, providing a clear quantitative hierarchy of circle similarity.

\begin{figure}[htbp]
    \centering
         \begin{subfigure}{0.32\textwidth}
         \centering         \includegraphics[width=\textwidth]{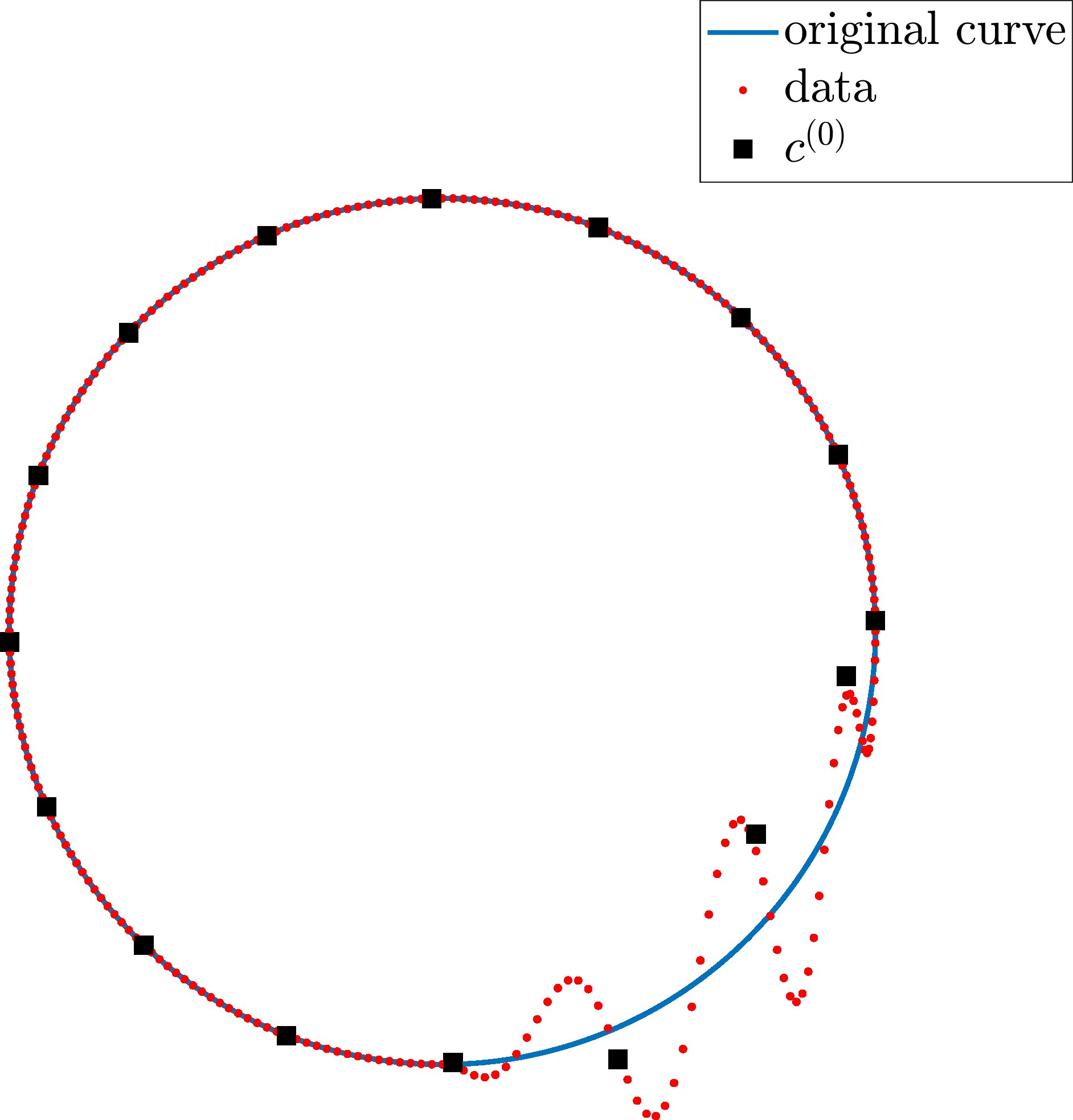}
    \end{subfigure}
    \qquad
    \begin{subfigure}{0.3\textwidth}
         \centering         \includegraphics[width=\textwidth]{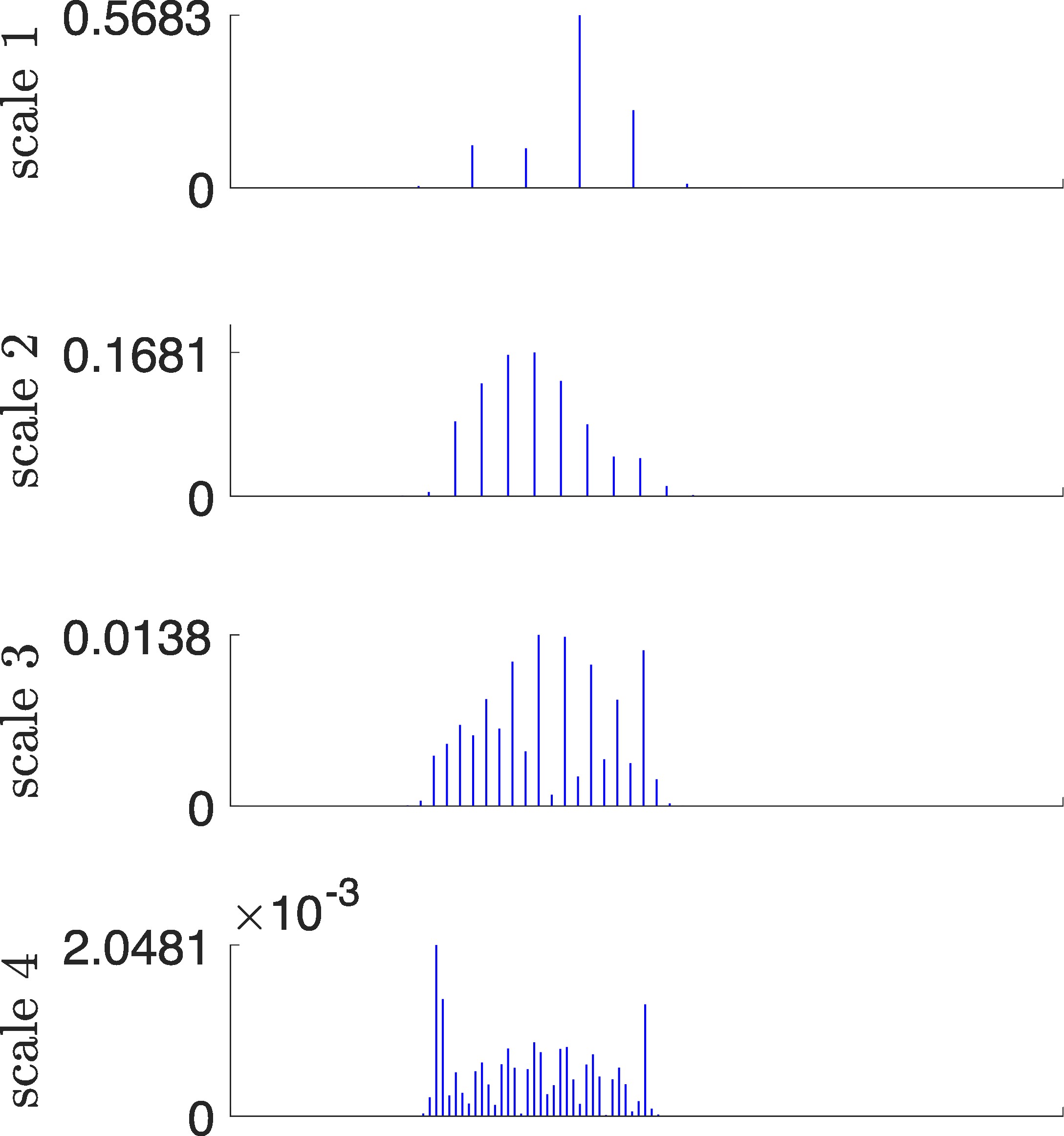}
     \end{subfigure}
     \caption{Anomaly detection in circles. On the left, $256$ equispaced data points are sampled from a circle with a lower right quadrant perturbed by a localized smooth function. Overlaid are 16 coarse data points (black squares). On the right, the multiscale representation via the Euclidean norms of the detail coefficients. It can be clearly observed that the magnitude of the detail coefficients increases dramatically exactly in the region affected by the perturbation.}
     \label{fig:quarter noisy circle pyramid}
\end{figure}

Since the multiscale transform relies fundamentally on local operators, it is inherently sensitive to localized variations and irregularities in the underlying data. To illustrate this sensitivity explicitly, in our second experiment we continue by demonstrating the efficacy of our geometric multiscale transform as an ``anomaly detector.''

Anomaly detection in geometric settings, particularly in planar shapes, involves identifying shapes or substructures that deviate significantly from a notion of normality within a population of 2-D geometries, which, in our example, are sampled circles. Unlike point-based anomaly detection, our approach operates at a higher level and aims to detect local structures that differ from the prescribed geometry. We expect the analysis to indicate such an anomaly in a local and consistent fashion across the analysis levels. 

In \figref{fig:quarter noisy circle pyramid}, we sample $256$ equispaced points from a perfect circle to which we add a smooth, localized oscillatory perturbation confined to a single quadrant. Upon applying the multiscale decomposition to this modified dataset, we observe a clear differentiation in the resulting detail coefficients. Specifically, the coefficients corresponding to the perturbed quadrant become substantially larger than those in the unperturbed regions, indicating that this region deviates from circular geometry, while the coefficients in the unperturbed areas remain close to numerical zero at double-precision accuracy. This distinction illustrates the transform's capability to detect, isolate, and characterize local irregularities, reflecting its intrinsic locality and spatial adaptivity.

\subsection{An illustrative analysis of network-generated samples}
\label{sec:NN detection}

Another application we consider is the detection and analysis of circles produced by a neural network. Generating exact geometric shapes is a known difficulty for learning-based generators~\cite{Chen2019, Park2019}: such models often produce shapes that appear nearly ideal, with deviations imperceptible to the naked eye, yet subtle discrepancies remain when examined from a mathematical perspective. Since the nonstationary pyramid transform introduced above provides a practical and sensitive tool for assessing the circularity of a given set of samples, we employ it as a proof of concept to distinguish between ``true'' circles, obtained by sampling a smooth underlying function, and approximate circles produced by a trained neural network.

\begin{figure}[htbp]
    \centering  
    \begin{subfigure}{0.35\textwidth}
         \centering         \includegraphics[width=\textwidth]{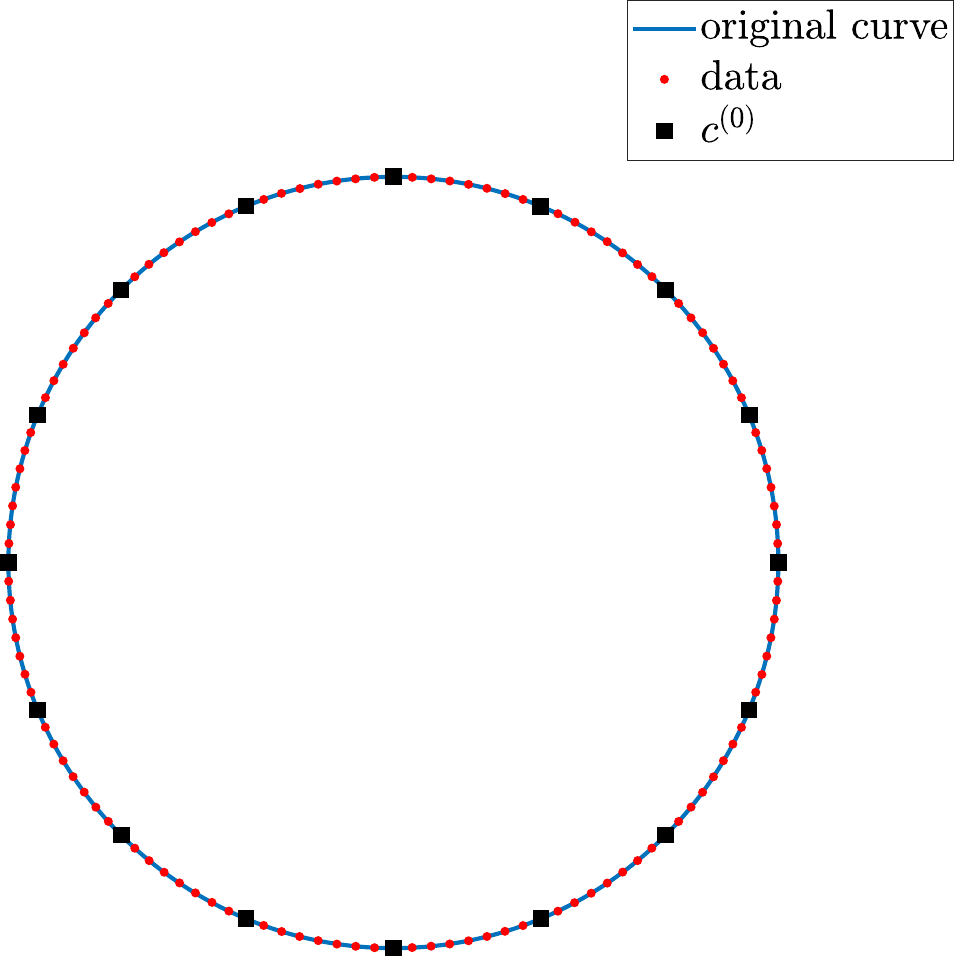}
    \end{subfigure}
    \qquad \qquad
    \begin{subfigure}{0.35\textwidth}
         \centering         \includegraphics[width=\textwidth]{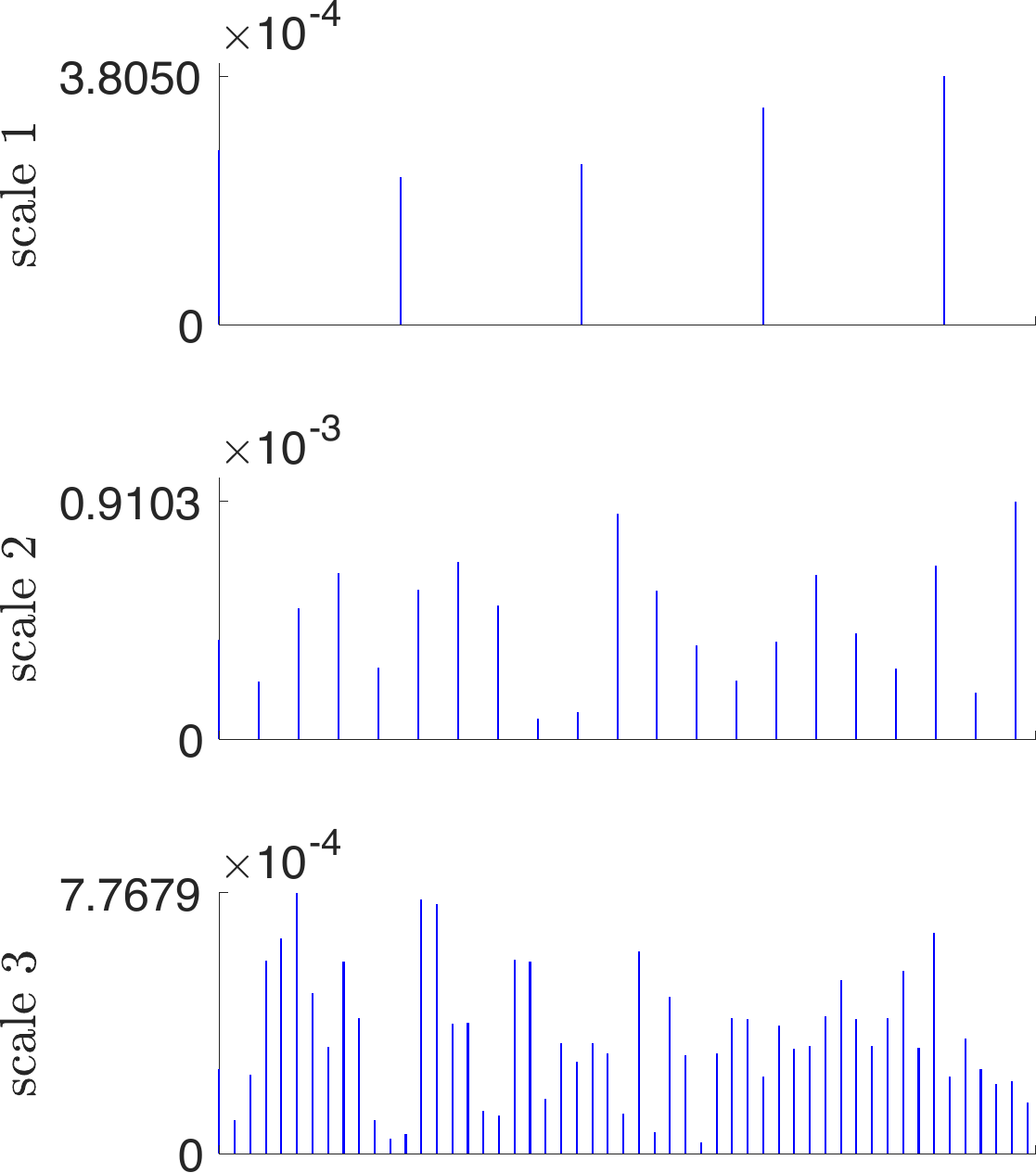}
     \end{subfigure}
    \caption{Multiscale representation of a circle generated by a neural network. On the left, $128$ points (red), obtained by providing a test center to a trained neural network and subsequently translating the output to be centered at the origin, are shown. The $16$ data points (black squares) depict the coarse approximation obtained by applying the analysis operator \eqref{NS multiscale transform} three times. On the right, the corresponding multiscale representation is displayed via the Euclidean norms of the detail coefficients. It is observed that the magnitudes of the detail coefficients are noticeably larger than those obtained from samples of an exact unit circle.}
    \label{fig:NN details}
\end{figure}

The circle prediction task was performed using a multi-layer feedforward neural network mapping a 2D center coordinate input to a $256$-dimensional output vector ($128$ points $\times 2$ coordinates). The network consists of three hidden layers with $32, 32$, and $16$ neurons, respectively. Training was conducted over $35$ epochs using an objective function constrained by a regularization parameter of $0.001$ to prevent overfitting. The dataset consisted of $10,000$ circles of radius $1$ with randomly selected centers in the domain $[0,10]\times[0,10]$. Each circle was uniformly sampled at $128$ points. The training was performed using the Levenberg--Marquardt algorithm~\cite{levenberg1944method, marquardt1963method}. After training, a test center was provided as input to the network, which then produced $128$ approximated sample points intended to lie on a circle of radius 1 centered at the given location. Subsequently, these points were translated by the prescribed center so that they are centered at the origin, allowing a direct comparison with the unit circle. The resulting samples, along with their corresponding detail coefficients obtained from the pyramid decomposition, are presented in \figref{fig:NN details}. 

Although these points appear to form a perfect circle from a visual perspective, the application of the proposed pyramid decomposition, based on the nonstationary scheme \eqref{NS geometric scheme}, reveals a markedly different behavior. In particular, the average detail norm $\nu_\ell(\boldsymbol{c}_{\text{NN}})$ is significantly larger across all levels $\ell$ than that obtained from exact circle samples $\nu_\ell(\boldsymbol{c}_{\text{clean}})$  (see \figref{fig:clean circle}). This indicates that the points do not lie precisely on the unit circle, but instead exhibit small yet measurable deviations from it.

The absence of monotone decay of $\nu_\ell(\boldsymbol{c}_{\text{NN}})$ over the available decomposition levels indicates that, relative to the chosen parametrization and subdivision model, the network error contains appreciable fine-scale oscillations. Since only finitely many samples and scales are considered, this behavior should not be interpreted as excluding the existence of a differentiable interpolating curve.

\begin{figure}[htbp]
    \centering
    \begin{subfigure}{0.22\textwidth}
         \centering         \includegraphics[width=\textwidth]{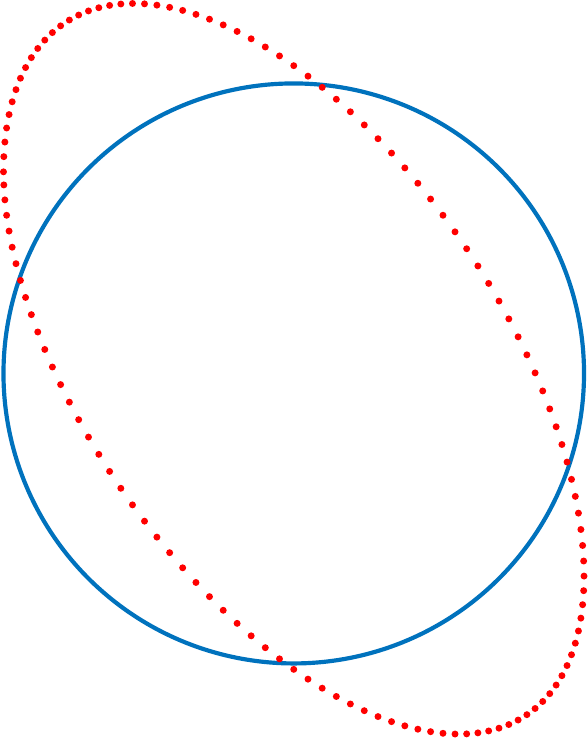}
     \end{subfigure}
     \hspace{2 cm}
     \begin{subfigure}{0.22\textwidth}
         \centering         \includegraphics[width=\textwidth]{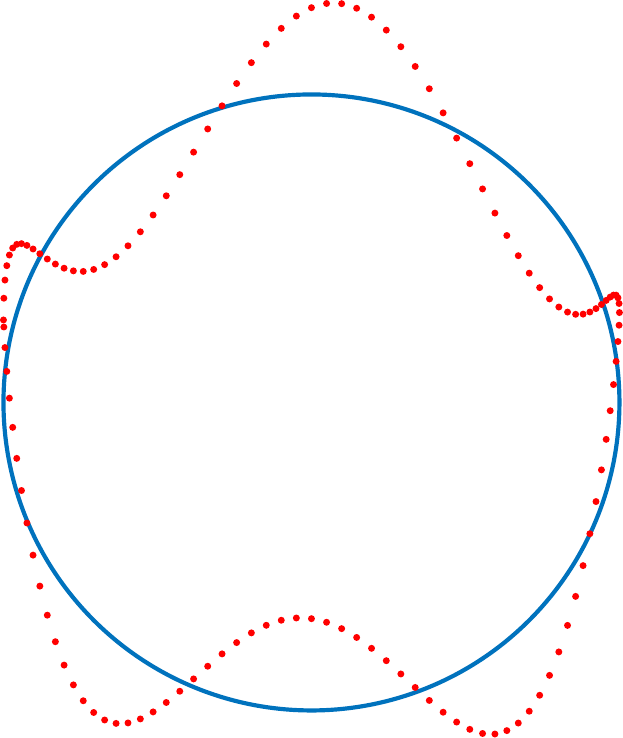}
    \end{subfigure}
    \hspace{2 cm}
    \begin{subfigure}{0.22\textwidth}
         \centering         \includegraphics[width=\textwidth]{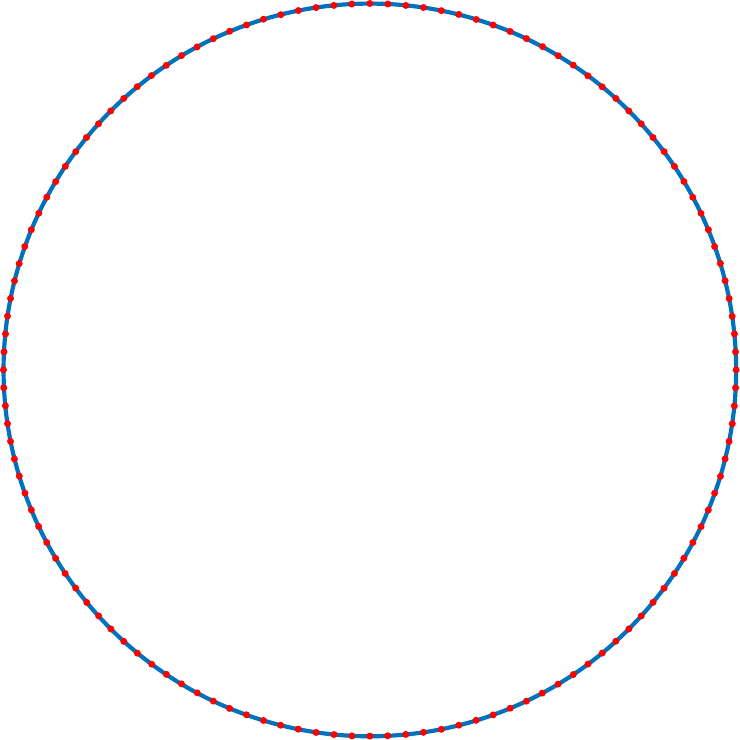}
    \end{subfigure} \\
     \begin{subfigure}{0.3\textwidth}
         \centering        \includegraphics[width=\textwidth]{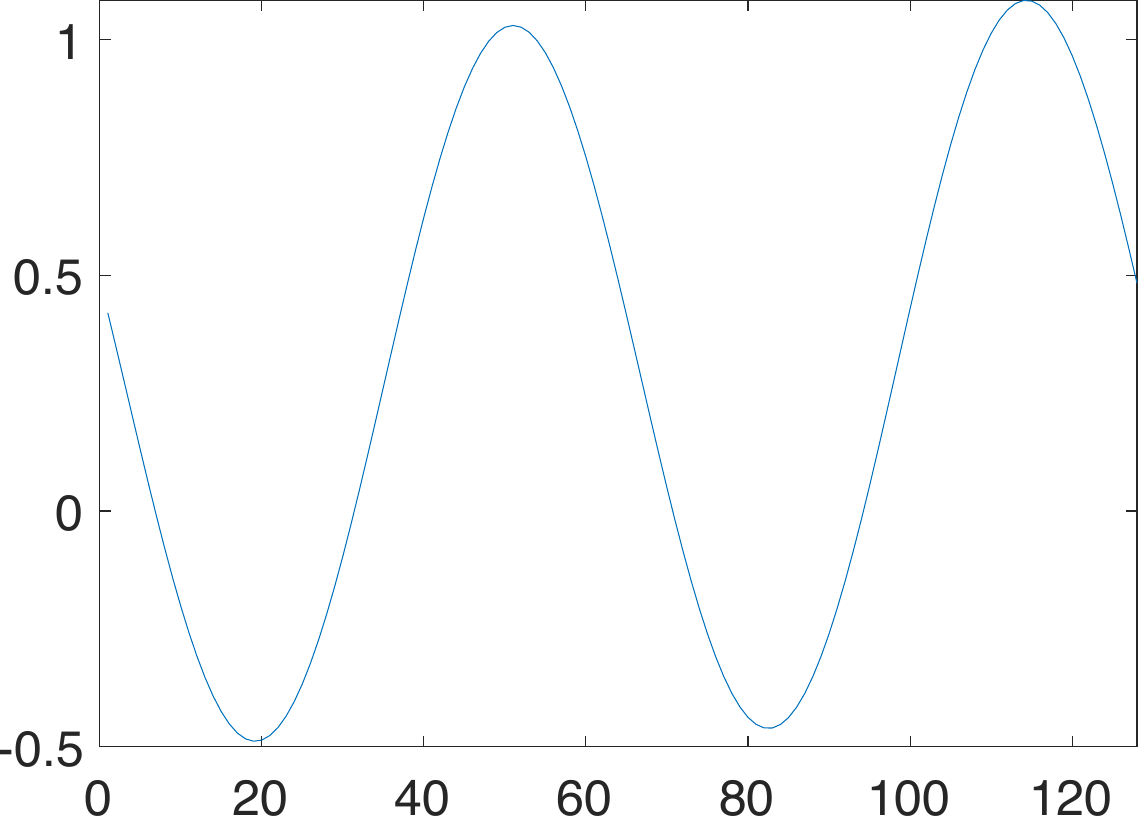}
    \end{subfigure}
    \hspace{0.5 cm}
    \begin{subfigure}{0.3\textwidth}
         \centering         \includegraphics[width=\textwidth]{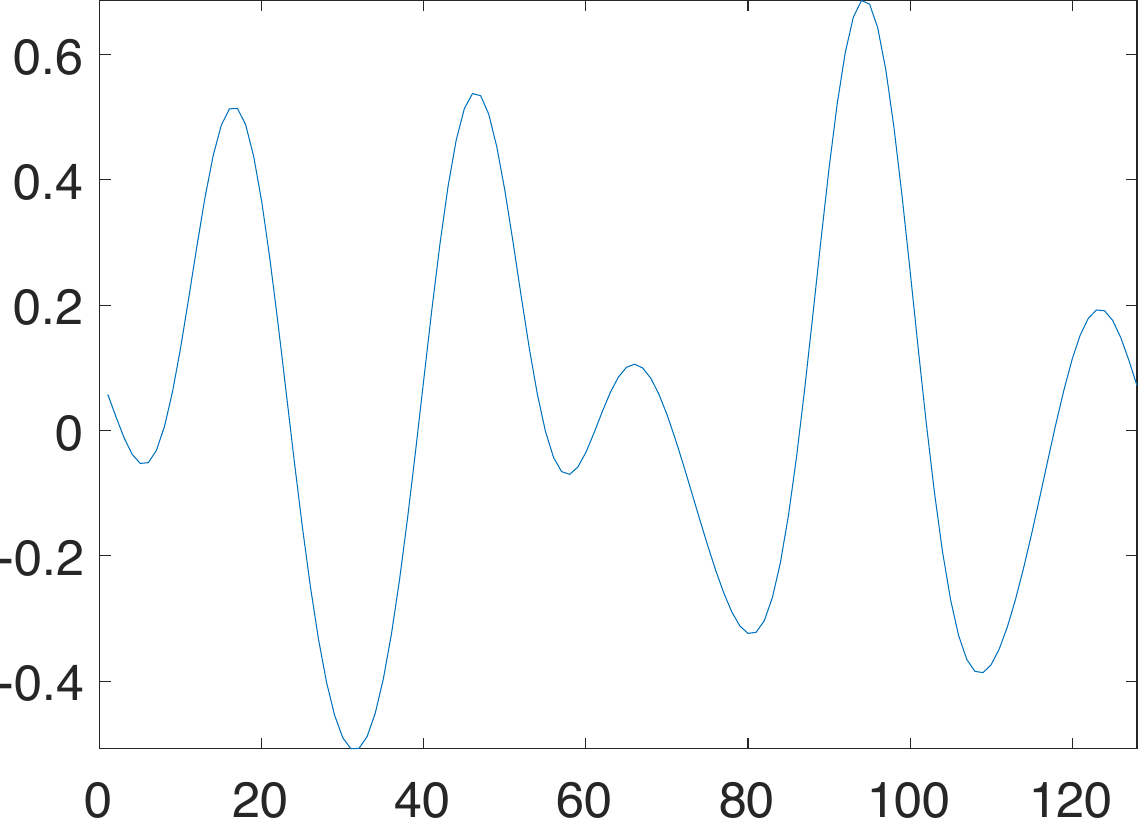}
     \end{subfigure}
     \hspace{0.5 cm}
    \begin{subfigure}{0.29\textwidth}
         \centering         \includegraphics[width=\textwidth]{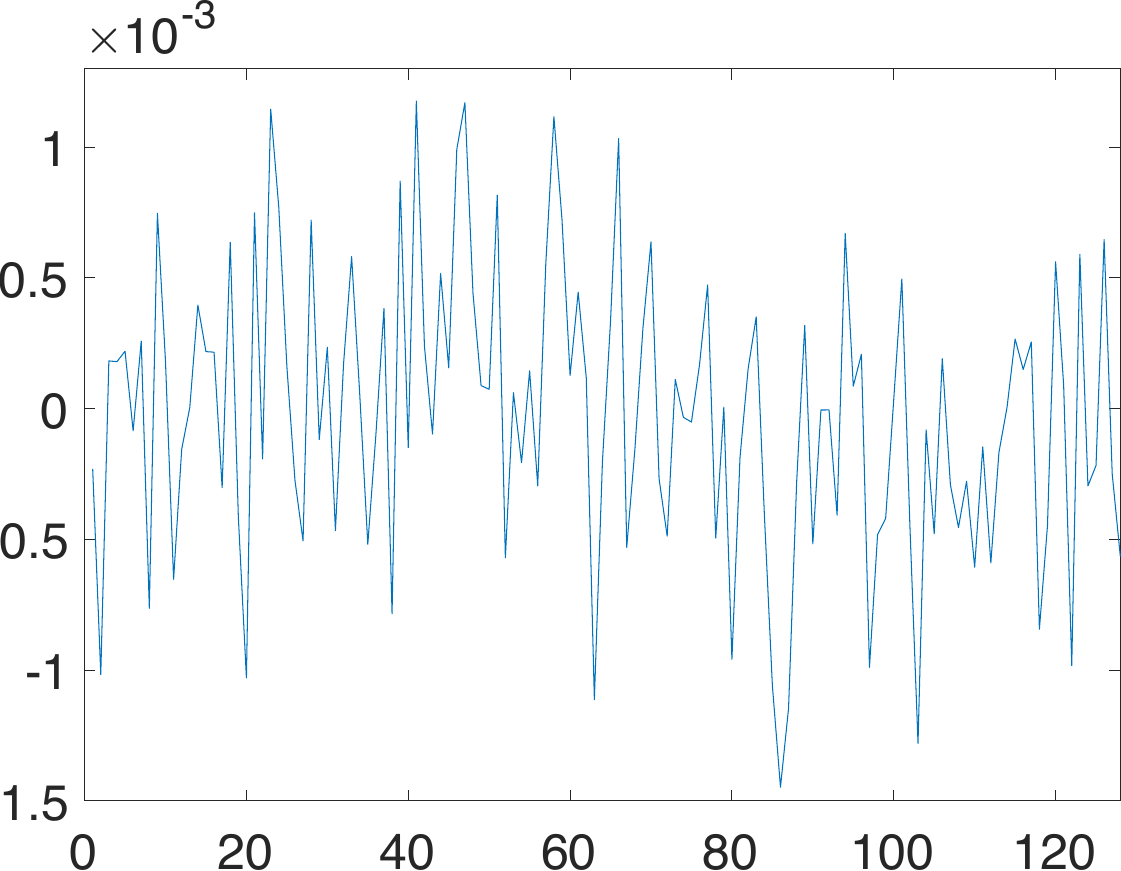}
     \end{subfigure}
     \caption{Quantification of pointwise deviations from the unit circle. The first two upper plots show $128$ points sampled from perturbed circular curves (``wavy circles''), while the third upper plot shows $128$ points obtained by the trained neural network shown in \figref{fig:NN details}. The lower plots illustrate the pointwise deviations of the sampled points from the unit circle. It is observed that points sampled from perturbed circular curves exhibit coherent and smooth deviations, whereas the points obtained via the trained neural network exhibit nonsmooth deviations.}
     \label{fig:diff from circle}
\end{figure}

To further illustrate this distinction, \figref{fig:diff from circle} presents the points obtained from the neural network learning process alongside two examples of points sampled from the unit circle and perturbed by smooth noise, together with their corresponding pointwise distances from the unit circle. It can be observed that the smoothly perturbed samples exhibit coherent and smooth deviations, whereas the neural-network-generated circle, despite displaying smaller absolute errors, is characterized by nonsmooth fluctuations.

Furthermore, as shown in \figref{fig:wavy_circle_pyramid} the smoothly perturbed circles yield substantially smaller values of $\nu_\ell(\boldsymbol{c}_{\text{smooth}})$ compared to $\nu_\ell(\boldsymbol{c}_{\text{NN}})$.
This observation demonstrates that $\nu_\ell(\boldsymbol{c})$ quantifies proximity to an ideal circular shape in terms of smoothness and deviation from the reproduced geometric model and the oscillatory character of that deviation, rather than mere visual similarity.

During training, the neural network generates samples that progressively approach a circular configuration as the number of epochs increases. This trend is explicitly reflected in the monotonic decrease of $\nu_\ell(\boldsymbol{c})$ across training epochs.

In our experiment, we targeted a maximal pointwise Euclidean error of $10^{-4}$. An accuracy of $10^{-3}$ was reached after $9$ epochs, whereas reducing the error from $10^{-3}$ to $10^{-4}$ required $26$ additional epochs. Thus, for the chosen architecture, regularization, optimizer, and stopping criterion, the rate of improvement slowed as training progressed. This behavior should not be interpreted as an intrinsic limitation of neural networks or of their number of trainable parameters. Indeed, because the sampling positions and radius are fixed in this experiment, the mapping from the prescribed center to the ordered circle samples is affine and can, in principle, be represented exactly by a suitable affine layer. The purpose of this experiment is therefore not to assess the expressive limitations of neural networks, but to demonstrate that the proposed multiscale transform can reveal small oscillatory artifacts produced by a particular trained model.

\begin{landscape}
\begin{figure} 
    \centering
    \begin{subfigure}{0.245\linewidth}
         \centering
         \includegraphics[width=\textwidth]{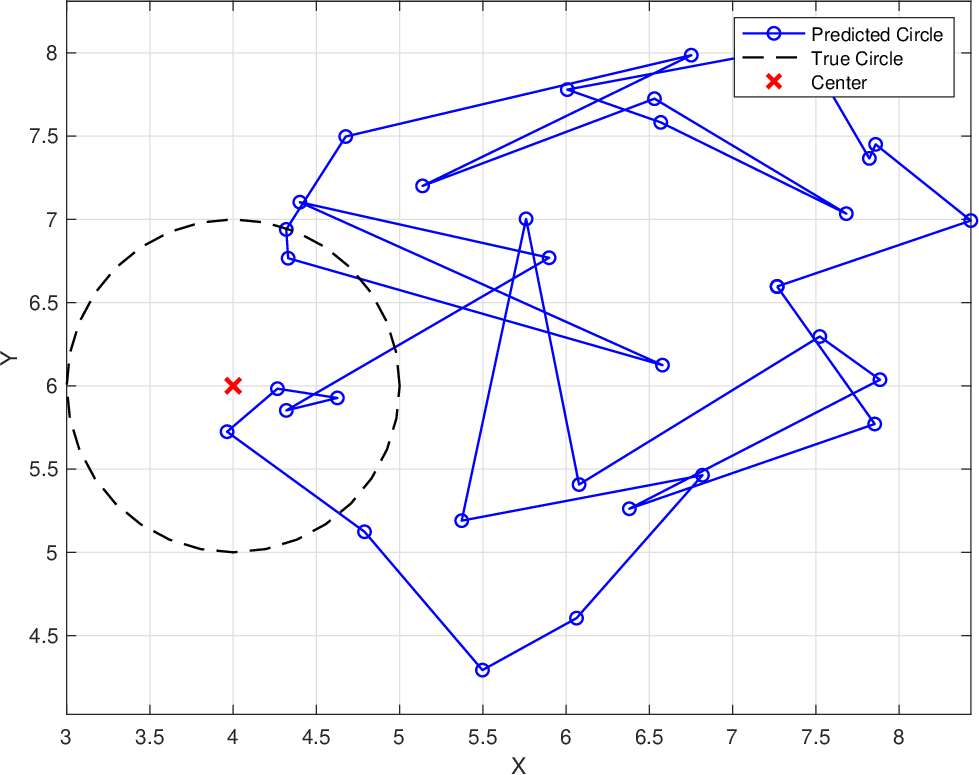}
    \end{subfigure}
    \begin{subfigure}{0.245\linewidth}
         \centering
         \includegraphics[width=\textwidth]{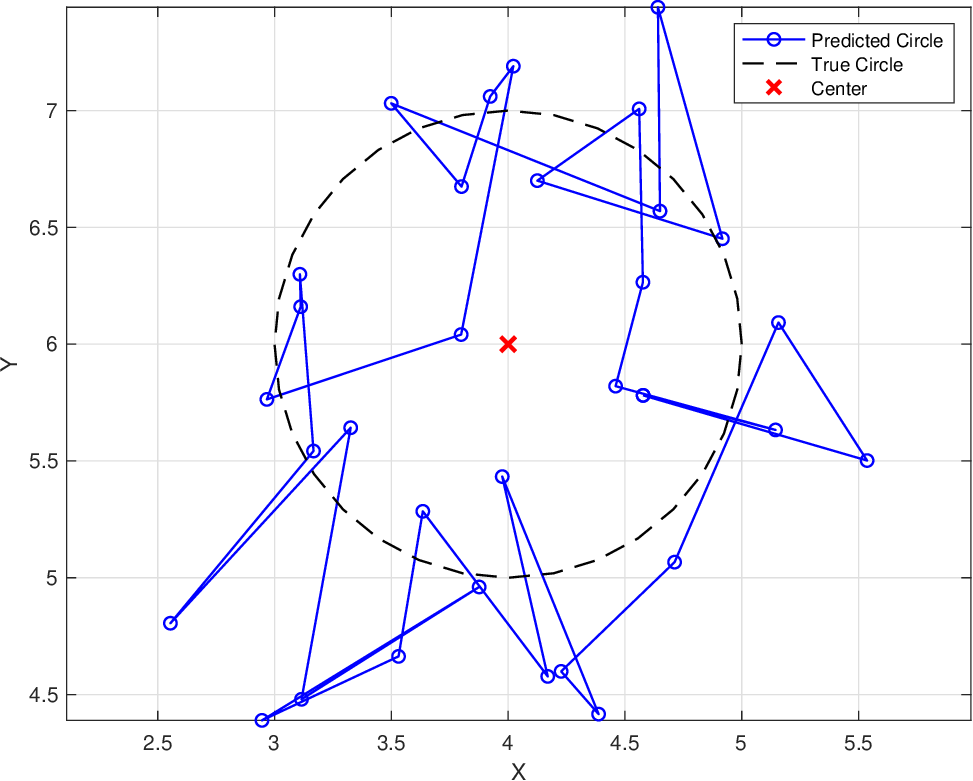}
    \end{subfigure}
    \begin{subfigure}{0.245\linewidth}
         \centering
         \includegraphics[width=\textwidth]{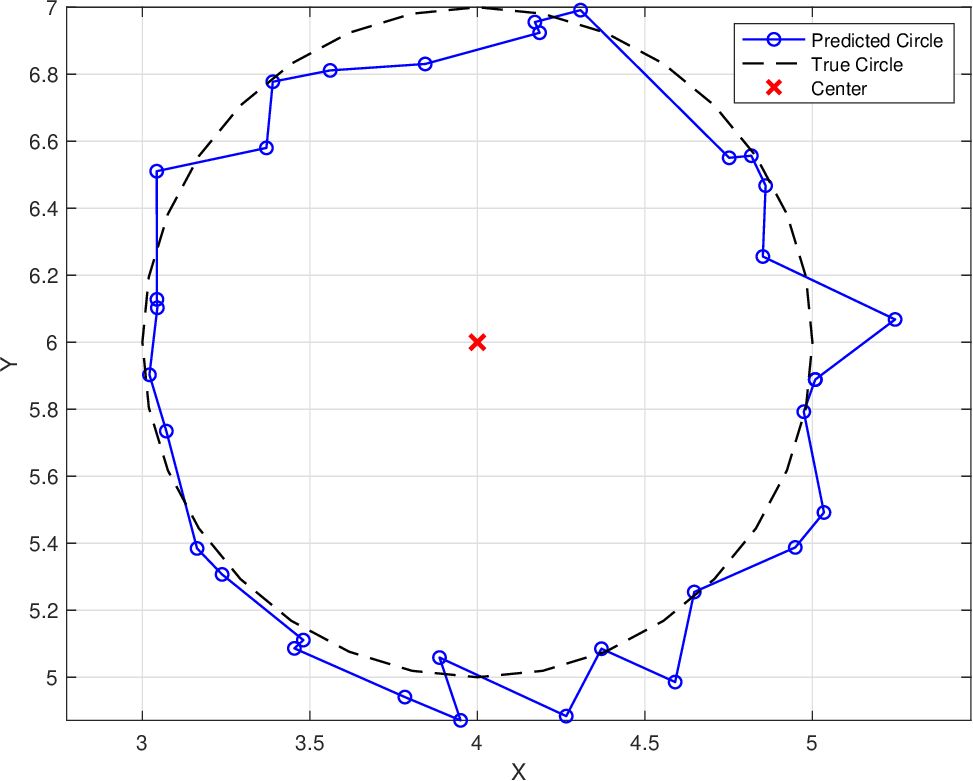}
    \end{subfigure}
    \begin{subfigure}{0.245\linewidth}
         \centering
         \includegraphics[width=\textwidth]{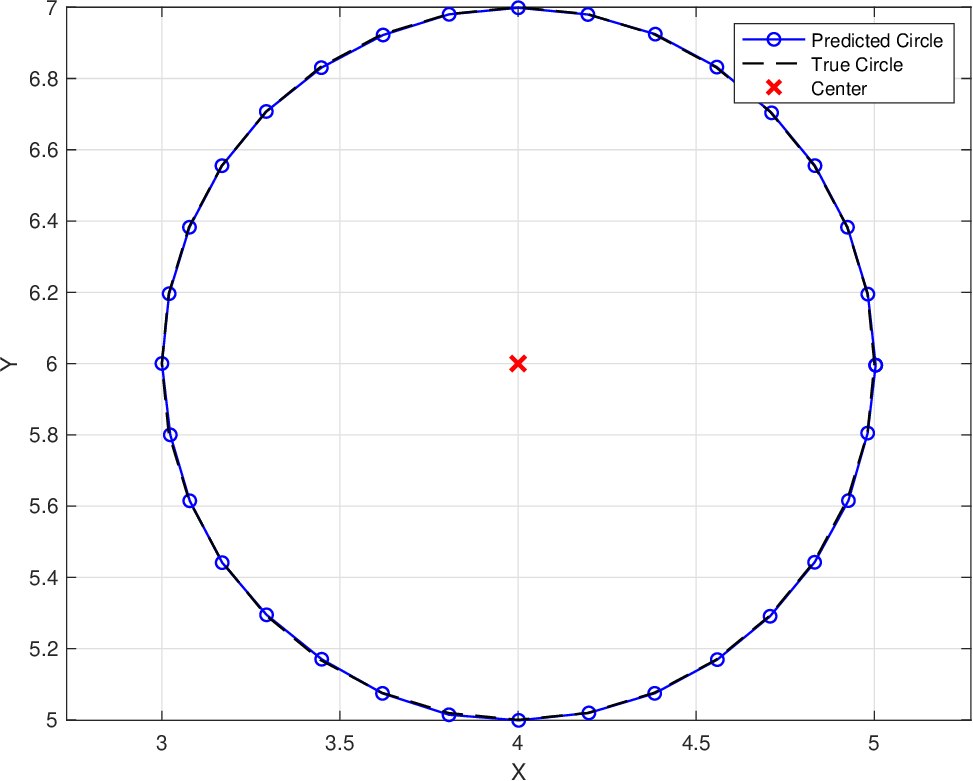}
    \end{subfigure}
   \\ \hfill
   \begin{subfigure}{0.2\linewidth}
         \centering
         \includegraphics[width=\textwidth]{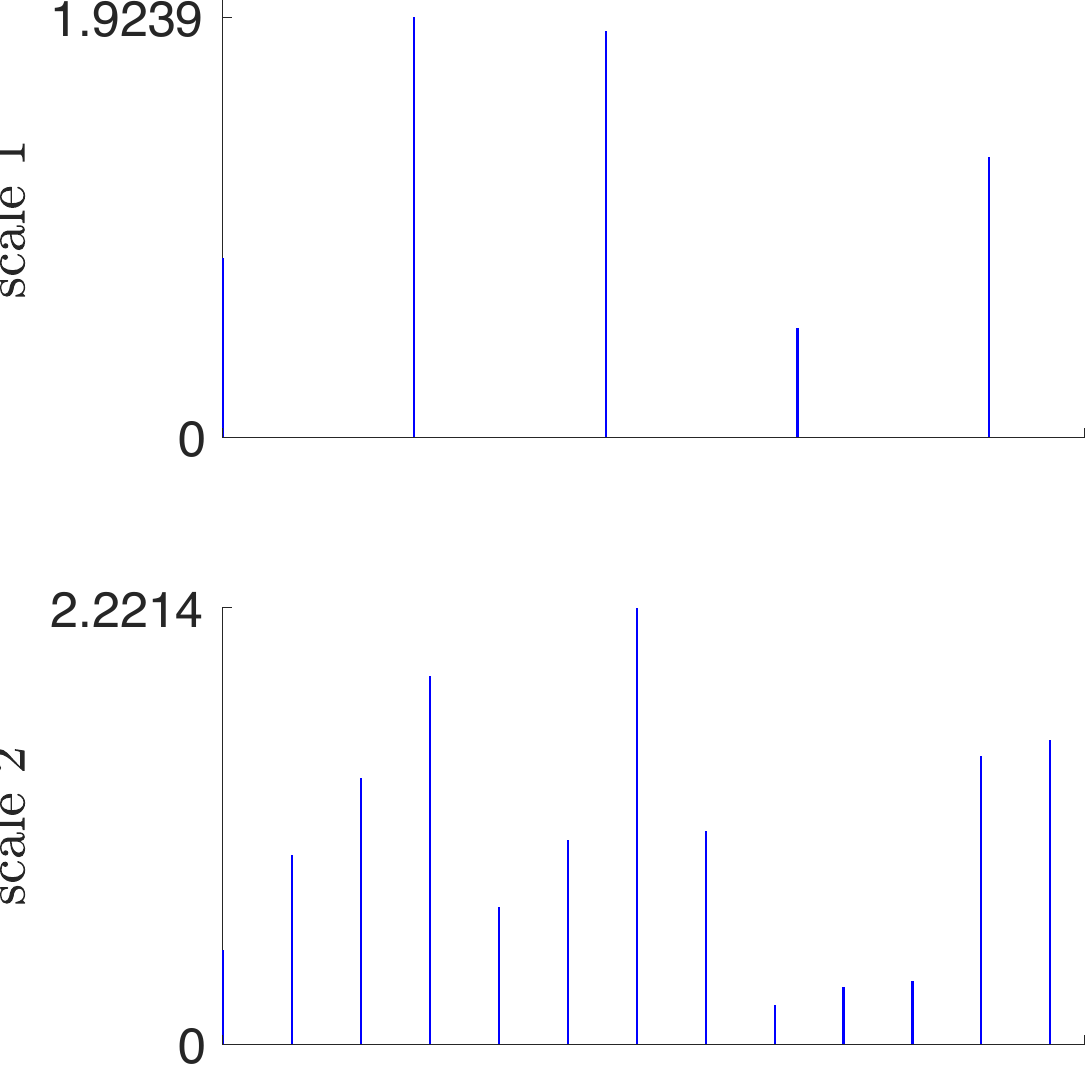}
         \caption{Epoch 1}
    \end{subfigure} \hfill
    \begin{subfigure}{0.2\linewidth}
         \centering
         \includegraphics[width=\textwidth]{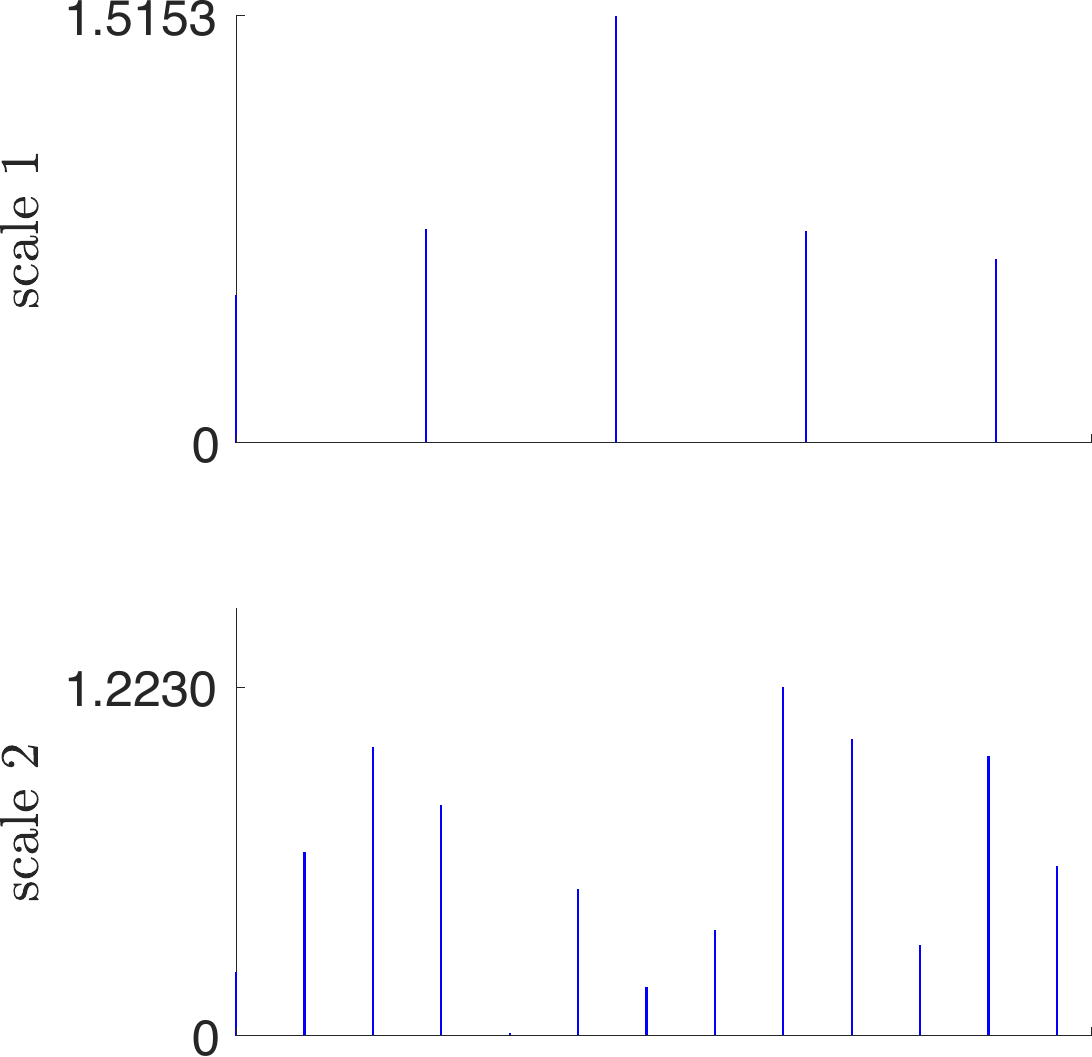}
         \caption{Epoch 3}
    \end{subfigure} \hfill
    \begin{subfigure}{0.2\linewidth}
         \centering
         \includegraphics[width=\textwidth]{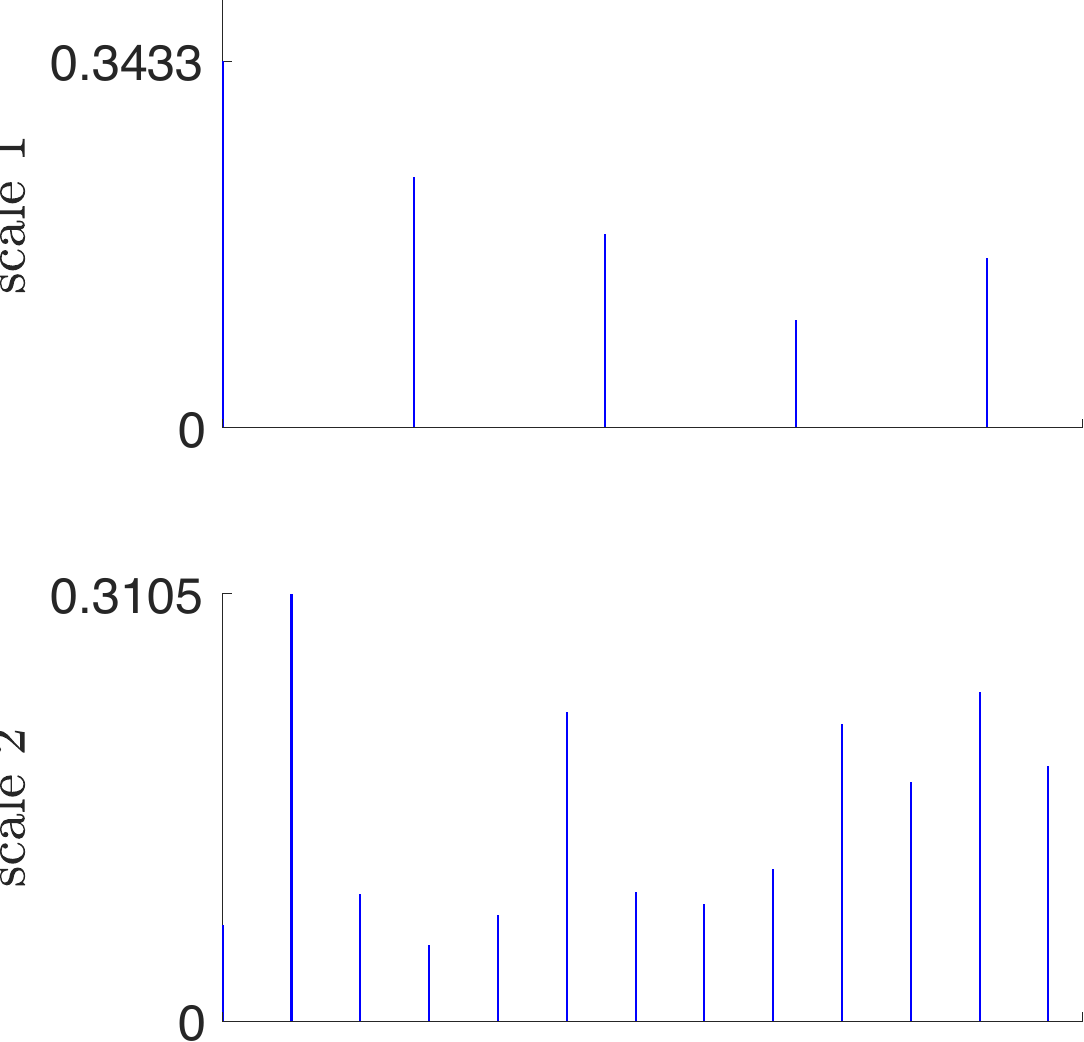}
         \caption{Epoch 5}
    \end{subfigure} \hfill
    \begin{subfigure}{0.2\linewidth}
         \centering
         \includegraphics[width=\textwidth]{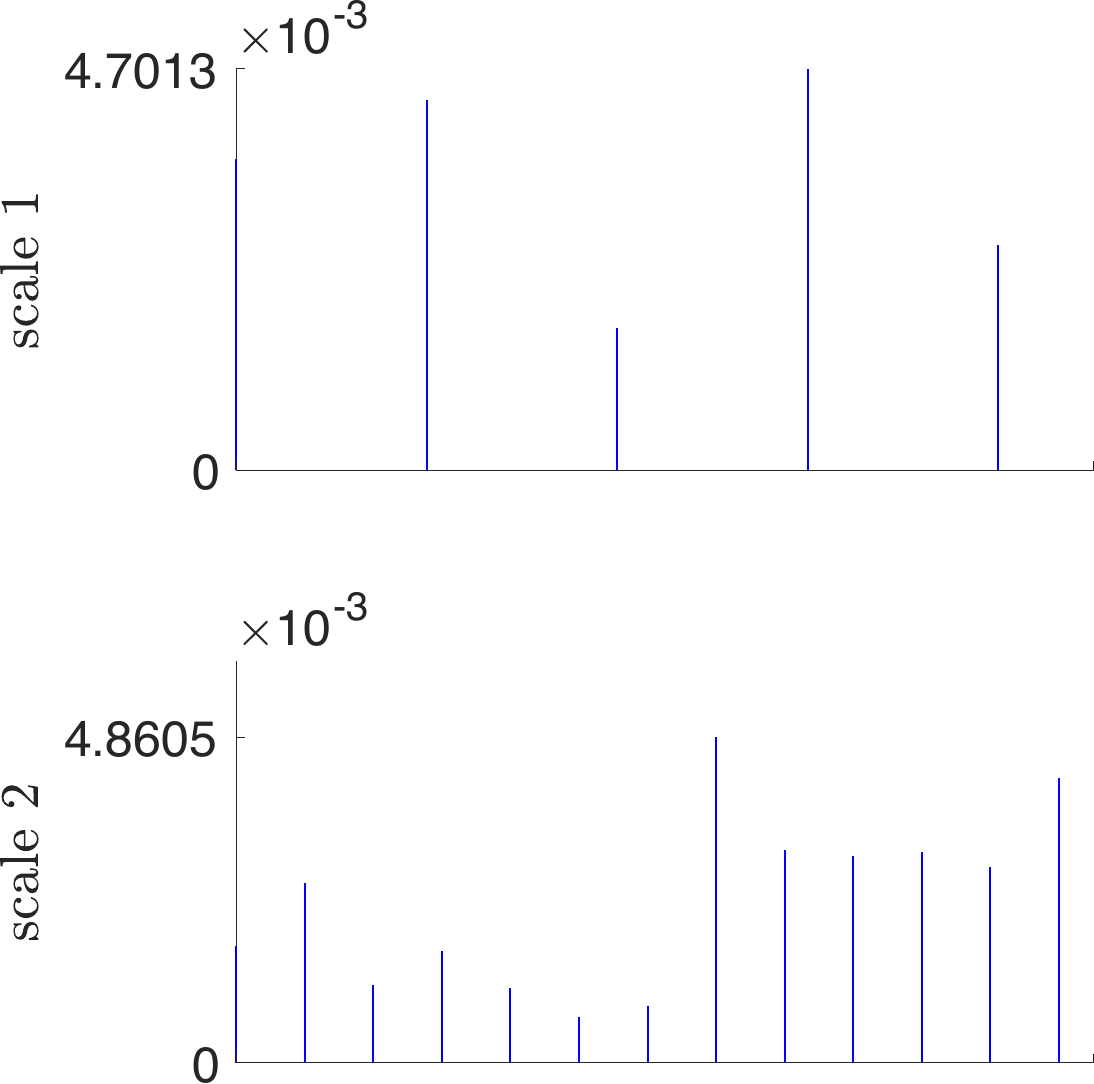}
         \caption{Epoch 7}
    \end{subfigure} \hfill
     \caption{Evolution of the predicted points and their corresponding detail coefficients throughout the training process. The predicted $32$-point approximations of a circle centered at $(4,6)$ after $1, 3, 5$, and $7$ training epochs are shown, together with their associated multiscale detail coefficients. In each of the panels (a)--(d), the predicted points are displayed on top and the Euclidean norms of the corresponding detail coefficients below. As the number of epochs increases, the values of $\nu_\ell(\boldsymbol{c})$ systematically decrease across all levels.}
     \label{fig:epochs}
\end{figure}
\end{landscape}

To analyze the network's learning dynamics, we conducted a follow-up experiment based on the previous configuration, with a key modification: the sampling density was reduced to $32$ points per circle. Due to this reduction, the training process converged and terminated after $7$ epochs. Rather than evaluating the network only after full convergence, training was periodically halted at epochs $1, 3, 5$, and $7$. At each checkpoint, the network was tested on a fixed center $(4, 6)$. The resulting $32$ approximated points were then processed by applying the pyramid transform \eqref{NS multiscale transform} based on subdivision scheme \eqref{NS geometric scheme} twice. The evolution of these generated points across the selected epochs is presented in \figref{fig:epochs}.

\section{Concluding remarks} \label{sec:conclusion}

In this paper, we have extended the theory of pyramid multiscale transforms by introducing nonstationary upsampling operators. Incorporating nonstationarity broadens the classical construction toward a more flexible class of pyramid transforms. To illustrate this, we presented detailed examples, proved fundamental properties, namely the decay of the detail coefficients and the stability of both the decomposition and the reconstruction, and demonstrated the practical benefits of the framework through applications focused on geometric features.

The applications highlight the role of locality, a property obtained directly from the upsampling operators that refinements of subdivision schemes. This locality both improves computational efficiency and lets the analysis adapt to varying geometric structure in the data, enabling a precise representation of geometric detail.

Our results indicate that the proposed framework has potential for a range of applications in geometric data analysis. While the examples presented here primarily serve as an initial demonstration, they illustrate the flexibility and practical value of the approach. Natural application areas include denoising, anomaly detection, feature extraction, image and shape analysis, and data compression, which makes the framework appealing for tasks involving geometric and structured data.

Looking ahead, several directions are worth exploring. These include extending the theory to nonlinear operators and adapting the approach to more complex, real-world data structures. Future work can also address applications involving manifolds, Lie groups, graphs, and higher-dimensional data. Finally, investigating computational improvements and integrations with modern machine learning techniques would further extend the relevance of this framework.


\section*{Acknowledgements}
NS and WM are partially supported by the DFG award 514588180. WM is partially supported by the Nehemia Levtzion Scholarship for Outstanding Doctoral Students from the Periphery (2023). NS is partially supported by the NSF-BSF award 2024791, the BSF award 2024266

\section*{Declaration of competing interest}
The authors declare that they have no known competing financial interests or personal relationships that could have appeared to influence the work reported in this paper.

\section*{Declaration of generative AI and AI-assisted technologies in the manuscript preparation process}
During the preparation of this work, the authors used ChatGPT (OpenAI) to assist with manuscript organization and language revision. After using this service, the authors reviewed and edited the content as needed and take full responsibility for the publication's content.

\section*{Code and data availability}
The code implementing the proposed transforms and reproducing all figures and examples of this paper is openly available at \url{https://github.com/HadarLandau/NSP}. No external datasets were used; all data analyzed in this work are generated by the accompanying code.

\bibliographystyle{elsarticle-num}
\bibliography{LSpyramid_bib}

\end{document}